\documentclass[sensors,article,accept,pdftex,moreauthors]{Definitions/mdpi} 
\usepackage[utf8]{inputenc}
\usepackage{textcomp}
\usepackage{chngcntr}
\usepackage{graphicx}
\usepackage{subcaption}
\usepackage{graphicx}
\usepackage{amsmath}
\usepackage[version=4]{mhchem}
\usepackage{siunitx}
\usepackage{longtable,tabularx}
\usepackage{amssymb}
\usepackage{nomencl}
\makenomenclature
\usepackage{etoolbox}
\renewcommand\nomgroup[1]{%
  \item[\bfseries
  \ifstrequal{#1}{A}{Sets, Indices, and Graph}{%
  \ifstrequal{#1}{B}{Parameters}{%
  \ifstrequal{#1}{C}{Decision Variables}{}}}%
]}
\usepackage{algorithm}
\usepackage{geometry}
\usepackage{atbegshi}
\usepackage{etoolbox}
\usepackage{xcolor}
\usepackage{soul}

\usepackage{booktabs} % for horizontal lines
\usepackage{array} % for adjusting column widths
\usepackage{caption}
\usepackage{longtable}

\usepackage{ifthen}

\usepackage{nomencl}
\usepackage{xpatch}

\usepackage{amssymb}

\newlength{\nomitemorigsep}
\firstpage{1} 
\pubvolume{1}
\issuenum{1}
\articlenumber{0}
\pubyear{2023}
\copyrightyear{2023}
\datereceived{ } 
\daterevised{ } % Comment out if no revised date
\dateaccepted{ } 
\datepublished{ } 
\hreflink{https://doi.org/} % If needed use \linebreak
\Title{Designing a Surveillance Sensor Network with Information Clearinghouse for Advanced Air Mobility}

\TitleCitation{Designing a Surveillance Sensor Network with Information Clearinghouse for Advanced Air Mobility}

\Author{Esrat F. Dulia $^{1}$ and Syed A. M. Shihab $^{2}$}

\AuthorNames{Firstname Lastname, Firstname Lastname and Firstname Lastname}

\AuthorCitation{Lastname, F.; Lastname, F.; Lastname, F.}
\address{%
$^{1}$ \quad Graduate Research Assistant, College of Aeronautics and Engineering, Kent State University, Kent, OH 44242, USA; edulia@kent.edu\\
$^{2}$ \quad Assistant Professor, College of Aeronautics and Engineering, Kent State University, Kent, OH 44242, USA; sshihab@kent.edu}

\corres{Correspondence: e-mail@e-mail.com; Tel.: (optional; include country code; if there are multiple corresponding authors, add author initials) +xx-xxxx-xxx-xxxx (F.L.)}

\abstract{To ensure safe, secure, and efficient advanced air mobility (AAM) operations, an AAM surveillance network is needed to detect and track AAM traffic. Additionally, a cloud-based surveillance data collection, monitoring, and distribution center is needed, where AAM operators and service suppliers, law enforcement agencies, correctional facilities and municipalities can subscribe to for receiving relevant AAM traffic data to plan and monitor AAM operations. In this work, we develop an optimization model to design a surveillance sensor network for AAM that minimizes total sensor cost while providing full coverage in the desired region of operation, considering terrain types of that region, terrain-based sensor detection probabilities, and meeting the minimum detection probability requirement. Moreover, we present a framework for low altitude surveillance information clearinghouse (LASIC), connected to the optimized AAM surveillance network for receiving live surveillance feed. Additionally, we conduct a cost-benefit analysis of the AAM surveillance network and LASIC to justify investment in it. We examine six potential types of AAM sensors and homogeneous and heterogeneous network types. Our analysis reveals the sensor types that are the most profitable options for detecting cooperative and non-cooperative aircraft. According to the findings, heterogeneous networks are more cost-effective than homogeneous sensor networks. Based on the sensitivity analysis, changes in parameters such as subscription fees, number of subscribers, sensor detection probabilities, and minimum required detection probability significantly impact the surveillance network design and cost benefit analysis.}


\begin{document}

\section{Introduction} \label{intro}

\subsection{Advanced Air Mobility}

AAM is envisioned to allow emerging short-haul aircraft, such as small uncrewed aircraft system (sUAS) and electric vertical takeoff and landing aircraft (eVTOL), to operate in the lower altitudes of the national airspace for passenger and cargo transportation and other use cases in the coming years. AAM is anticipated to offer a number of benefits to society and environment over traditional ground transportation systems, including a considerable reduction in travel and delivery times, increased operational safety, and reduced negative impact on the environment \cite{rothfeld2021potential, dulia2021benefits}. Federal agencies, such as NASA and FAA, and industry and academia have been focusing their research on AAM aircraft design \cite{silva2018vtol}, concepts of operation \cite{thipphavong2018urban}, air traffic management \cite{FAA2020}, trajectory planning \cite{pradeep2019energy}, deconfliction \cite{yang2018autonomous}, market studies \cite{hasan2018urban, reiche2018urban}, network planning \cite{german2018cargo, lim2019selection, chen2022scalable} and operations planning \cite{shihab2019schedule, shihab2020optimal}. A more comprehensive review of past and recent AAM research can be found in \cite{garrow2021urban, straubinger2020overview}.

\subsection{Motivation and Contributions} 
% or Research Needs
\subsubsection{Surveillance Sensor Network Design for Advanced Air Mobility}

\begin{figure}[H]
   
\begin{adjustwidth}{-\extralength}{0cm}
%\centering %% If there is a figure in wide page, please release command \centering
 \centering
    \includegraphics[width=12cm,height=8.5cm]{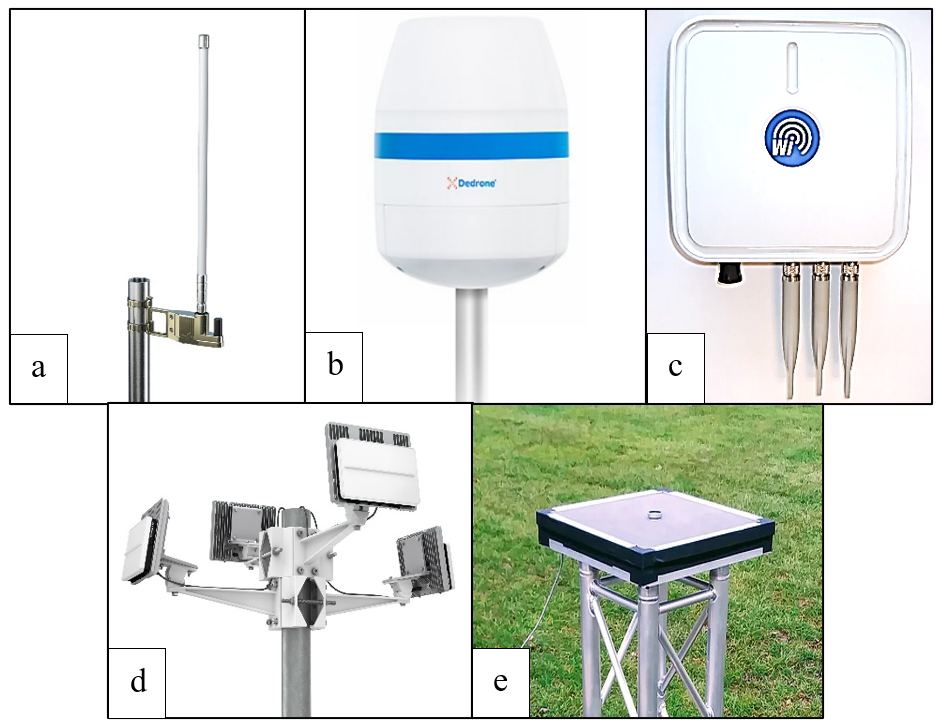}
\end{adjustwidth}
    \caption{Sensors of different types: (a)  ADS-B receiver \cite{ads-b}, (b) radio frequency sensor \cite{RF360}, (c) remote ID receiver \cite{remote}, (d) radar \cite{echo}, and (e) acoustic sensor \cite{aco}.}
    \label{fig:comb}
\end{figure}

While AAM research has been advancing on many fronts, one area of research which is critical to enabling AAM but has not received much attention is surveillance sensor network design for AAM. New surveillance sensor networks are needed specifically for AAM to detect and track AAM traffic to ensure efficient, safe and secured AAM operations. Much of the existing surveillance infrastructure for conventional aviation is not adequate for AAM for mainly two reasons. Firstly, AAM is envisioned to involve operations of aircraft within new urban and suburban operating regions, where no sensors currently exist for aircraft surveillance. Secondly, traditional aviation sensors which may already exist in anticipated operating regions for AAM will not be adequate because they have not been designed specifically to detect small aircraft, such as sUAS and eVTOLs, and accurately identify multiple aircraft flying near each other at lower altitudes in inclement weather conditions \cite{FAAdronedetection}, as would be the case for future high density AAM operations. Hence, AAM-specific sensors with specialized features and capabilities are needed that enable real-time detection and tracking of AAM aircraft in various weather conditions, simultaneous detection of multiple aircraft, and accurate identification and classification of AAM aircraft. Keeping these requirements in mind, several different types of surveillance sensors for AAM --- including, radar, radio frequency sensor, Automatic Dependent Surveillance-Broadcast (ADS-B) sensor, remote ID sensor, optical sensor, and acoustic sensor --- have been developed by various sensor manufacturers, such as Echodyne, Dedrone, and AVIONIX. Examples of such sensors are pictured in Figure \ref{fig:comb}. These sensors use either electromagnetic or sound waves to determine the distance, angle, and radial velocity of aircraft relative to their installation sites to detect and track aircraft. Such AAM sensors need to be set up to form a surveillance network in future AAM operating regions. While other research has focused on designing surveillance networks for both aviation and non-aviation applications, the specific context of AAM remains largely unexplored. To bridge this gap, we have developed a Surveillance for AAM Network Design (SAND) optimization model for identifying the optimal locations for placing the sensors to build the AAM surveillance network such that: 1) full coverage is provided in the desired region of operation; 2) minimum detection probability requirement is satisfied; and 3) the total sensor cost is minimized. In determining the optimal sensor placement solution, the model considers the range of various sensor types and the degradation of detection probabilities based on terrain types within the operating region. We consider Ohio as a case study for designing an AAM surveillance network using this model.

% Six sensors of different types: (a) EchoGuard Radar \cite{echo}, (b) PingStation ADS-B receiver \cite{ads-b}, (c) Drone Scout remote ID receiver \cite{remote}, (d) Dedrone sensor RF-360 \cite{RF360}, (e)  Drone Hound acoustic sensor \cite{aco}, and (f) Q6225-LE PTZ Network camera \cite{camera}

\subsubsection{Low Altitude Surveillance Information Clearinghouse}

The demand for AAM is expected to grow rapidly in the coming years due to several factors such as urbanization, population growth, and the ever increasing need for more efficient and sustainable transportation solutions \cite{goyal2021advanced}. Therefore, significant amounts of AAM traffic surveillance data will be generated by the AAM surveillance network, which would require efficient storage and processing solutions to ensure that the data is easily accessible and available for real-time and offline analysis by relevant AAM stakeholders, such as AAM operators, airspace service providers and law enforcement agencies. A digital LASIC can act as a central repository for this traffic data, allowing for data accessibility and sharing among various entities for flight planning, aircraft routing, air traffic control, counter uncrewed aircraft systems (UAS) operations planning, and better coordination among low altitude airspace users. Some of the functions LASIC can enable for its users include: access to live surveillance feeds, real-time coverage map, and archival data; data analytics and visualization; tactical deconfliction; and querying current and historical UAS positions by UAS ID and by location. 

For implementing and hosting LASIC, a cloud server is considered to be more suitable than a local server as cloud computing can provide a scalable, flexible, and cost-effective platform for ingesting, processing, storing, analyzing, and sharing large amounts of transportation data \cite{shengdong2019intelligent} generated by AAM traffic. Cloud computing can improve the performance and efficiency of transportation systems such as LASIC by relocating the hardware and software components to the cloud network \cite{nayar2018cost}, which would allow LASIC to access the computing resources and data storage capabilities of the cloud network. This can potentially reduce the need for expensive hardware and infrastructure associated with local servers on-site, while also providing greater flexibility and scalability for LASIC operations.

% Hence, cloud computing is suitable for handling the anticipated big data associated with AAM traffic \cite{shengdong2019intelligent}. 

% With a focus on efficient and scalable use of computing resources and commitment to environment-friendly information technology, cloud computing is quickly becoming a prominent paradigm for data storage and computation. 

\begin{figure}[H]
\begin{adjustwidth}{-\extralength}{0cm}
\centering %% If there is a figure in wide page, please release command \centering
\includegraphics[width=15cm,height=8cm]{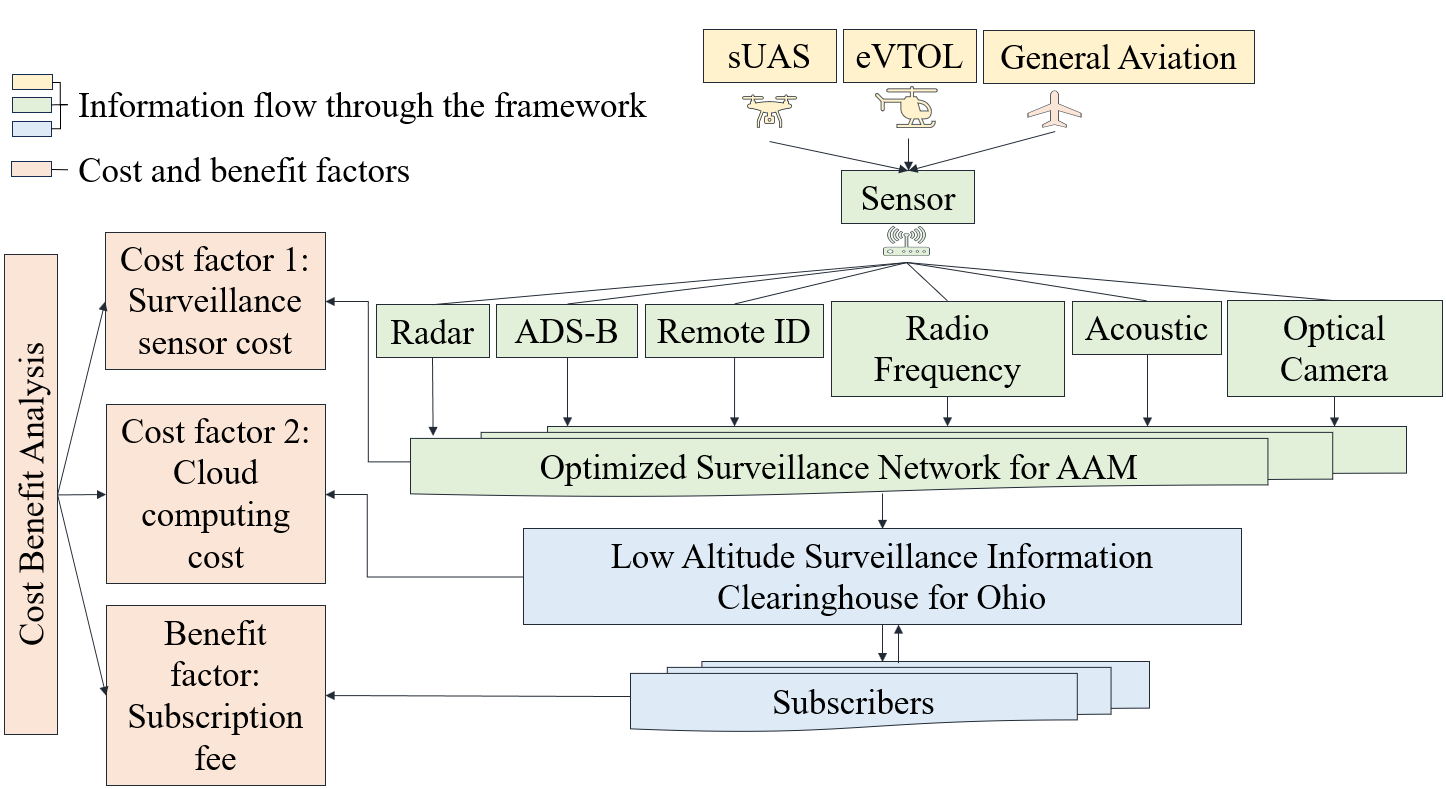}
\end{adjustwidth}
    \caption{Overview of LASIC framework and associated cost and benefit factors.}
    \label{fig1}
\end{figure}

An overview of the optimized AAM surveillance network and LASIC framework and its associated cost and benefit factors are illustrated in Figure \ref{fig1}. Based on a survey of the present AAM sensor market, we selected six different sensor types: radar, radio frequency sensor, ADS-B, remote ID, optical camera, and acoustic sensor. The surveillance and telemetry data associated with sUAS, eVTOL, and general aviation traffic --- such as position, velocity, flight intent, remote identification --- can be captured and generated by the optimized surveillance network, allowing the aircraft movement in the airspace to be tracked. This surveillance data can then be ingested into LASIC, which will provide the subscribers of LASIC with information about scheduled and real time AAM operations and relevant airspace activities so that they may plan for their flight operations accordingly. The subscribers of LASIC will potentially include: AAM operators engaged in different AAM use cases such as passenger and cargo transportation, bridge inspections, medical delivery, etc.; airspace service providers; law enforcement agencies; correctional facilities; and municipalities.

As for any other major transportation infrastructure project, to justify the investment in AAM surveillance network and LASIC, a rigorous cost-benefit analysis is needed \cite{couture2016cost, Mishan2020, uckelmann2012performance}. Such an analysis is crucial to identify, quantify, and evaluate costs and benefits associated with the surveillance network and LASIC. We conduct this cost-benefit analysis for the State of Ohio by analyzing the associated cost and benefit factors of LASIC. The analysis period was considered to be the next 10 years, from 2024-2033. The two major cost factors of LASIC considered are: 1) surveillance sensor cost, the cost to purchase the sensors needed in the AAM operating regions in Ohio, which we estimate based on the results generated from SAND model; and 2) cloud computing cost to store and process the surveillance data in LASIC. The monthly subscription fee that a subscriber will pay to get access to use LASIC features is considered as the main benefit factor in this analysis. The cost-benefit analysis can be used to estimate the break-even point (BEP) for the different sensor types, the time to reach break-even in terms of net present value (NPV) of the return generated in the AAM operating regions.

% \subsection{Summary of Contribution}

\subsubsection{Summary of Contributions}

This paper addresses the critical need for surveillance network design in the rapidly emerging landscape of AAM, offering insights into sensor selection, network optimization, data management, and economic feasibility. Our key contributions are as follows:

\begin{enumerate}%[label=\alph*)]

% \item \hl{We underscore the need for novel surveillance networks exclusively designed for AAM. Conventional aviation surveillance infrastructure falls short due to the distinctive urban and suburban settings, coupled with the presence of AAM aircraft like sUAS and eVTOLs.}

\item[(a)] We develop the SAND model, which can determine optimal locations for sensor deployment to design a comprehensive AAM surveillance network that minimizes total sensor cost. The SAND model can provide full coverage in the desired AAM operating regions, considers terrain types within those regions, terrain-based sensor detection probabilities, and minimum detection probability requirements. We consider the State of Ohio as our case study and apply the SAND model to design an AAM surveillance network there.
  
  % sensor coverage and data requirements, sensor detection probability on different terrains, and sensor costs.}

 % \item \hl{We consider several sensor types tailored to AAM prerequisites acknowledging the limitations of existing sensors. These encompass radar, radio frequency, ADS-B, remote ID, optical, and acoustic sensors.}
 
\item[(b)] We consider several sensor types, such as radar, radio frequency, ADS-B, remote ID, optical, and acoustic sensors, to design two types of AAM surveillance sensor networks: homogeneous and heterogeneous. Our analysis of homogeneous sensor placement indicates that ADS-B and remote identification sensor types are the most profitable options for detecting cooperative aircraft, whereas the radio frequency sensor type is the most profitable option for tracking both cooperative and non-cooperative aircraft. According to the findings, implementing a heterogeneous sensor network composed of various sensor types is more cost-effective in reducing the overall sensor cost compared to a homogeneous sensor network that employs only one type of sensor.

\item[(c)] We present a cloud-hosted LASIC framework, which allows for managing and sharing of AAM surveillance traffic data. We compute the cost of operating the framework considering AAM traffic projections and relevant surveillance data generated in the AAM operating regions in Ohio, surveillance data types, interface standards, data sizes, cloud components, and cloud
  pricing policies. 

\item[(d)] We conduct a rigorous cost-benefit analysis of the proposed AAM surveillance network and LASIC implementation for the State of Ohio to determine the break-even points for different sensor types. We consider the uncertainty associated with AAM demand to determine the possible range of of costs, revenue, and NPV for the AAM surveillance network and LASIC.

\item[(e)] We perform a sensitivity analysis on key parameters of our study, including subscription fees, number of initial subscribers, terrain-based sensor detection probabilities, and minimum required detection probability. The insights demonstrate that changes in these parameters significantly impact the number of sensors required, total sensor cost, and  NPVs of results generated from the study. Our study provides policymakers with valuable insights to make informed decisions regarding investment in an AAM surveillance network and LASIC.

% \item \hl{The insights from our study demonstrate that higher terrain-based sensor detection probabilities result in a decrease in the number of sensors and the overall sensor cost. Additionally, an increase in the minimum required detection probability leads to a increase in the number of sensors and an associated increase in total sensor cost. Among sensor types, the acoustic sensors are most impacted by this increase, while the ADS-B sensors are the least affected. Furthermore, the analysis indicates that higher subscription fees and a larger number of subscribers contribute to higher NPV generated by LASIC.}

\end{enumerate}

\subsection{Outline of the Paper}

The remainder of this article is structured as follows. The relevant literature is analyzed in Section \ref{sec2}. In Section \ref{sec3}, the SAND model is presented and the potential cost and benefit factors of LASIC are discussed. After that, the results are presented and analyzed in Section \ref{sec4}. The paper is finally concluded in Section \ref{sec5} with the summary of insights gained from the analysis and potential extensions of this study.

\section{Literature Review}\label{sec2}

%revised

% and prospective directions for future AAM research
% For example
% addressed

% \hl{Several research efforts have been dedicated towards the recent development of AAM. In } \cite{straubinger2020overview}, \hl{the authors conducted a comprehensive review of the current state of AAM development. The authors discussed unresolved challenges related to regulations, infrastructure, and economics by synthesizing insights from various segments of the AAM research community. Their review encompassed diverse aspects, including aircraft specifications, integration challenges, certification, policy, traffic management, ground infrastructure, operational strategies, market structures, integration with existing transportation systems, and public acceptance. Additionally,} \cite{garrow2021urban} \hl{performed a meta-analysis of approximately 800 articles published between January 2015 and June 2020 in the domains of air, electric, and autonomous vehicles, comparing research themes to guide future AAM research directions. The authors also examined articles focusing on demand modeling, operational considerations, and the integration of AAM with existing urban infrastructure, utilizing insights from both the meta-analysis and comprehensive review to shape the identification of forthcoming AAM research priorities.}

A review of research on AAM surveillance and the general surveillance network design problem is presented in this section. For a more broader overview of AAM research, interested readers are referred to \cite{straubinger2020overview} and \cite{garrow2021urban}, where the authors collectively discussed prior AAM research and unresolved AAM challenges related to aircraft specifications, regulations, certification, policy, demand modeling, traffic management, ground infrastructure, operational strategies, market structures, integration with existing transportation systems, and public acceptance.

\subsection{AAM Surveillance}

A number of recent studies in the literature have concentrated on surveillance technologies, frameworks, and simulations aimed at tracking and monitoring AAM aircraft. Notably, NASA's recent work (\cite{kawamuraground}) has emphasized the necessity for surveillance of AAM aircraft, underlining the difficulty of modifying present air traffic control and management systems to accommodate the increased number of AAM aircraft in the lower airspace. In response, they has introduced a ground-based vision tracker that employs a vision tracking method with fixed cameras to monitor airborne objects, effectively sidestepping issues related to electromagnetic interference. Additionally, \cite{kannan2023simulation} developed a simulation system to model and assess AAM flight operations in densely populated urban areas using both air and ground-based sensors such as radar, LiDAR, and vision-based sensors. The purpose of this paper was to present the architecture and simulation setup for evaluating airborne autonomy technologies for urban AAM operations. In \cite{belwafi2022unmanned}, the authors mentioned remote identification as an emerging technology that allows ground observers to identify drones within airspace. Their objective was to provide a comprehensive overview and tutorial of the current status of regulatory, standardization, design, implementation, and testing efforts in the field of remote identification technology. In another study (\cite{lofu2022uranus}), the authors developed a surveillance framework to address the growing security threats facing critical infrastructures such as airports, military bases, city centers, and other restricted zones. This framework utilizes radio frequency (RF) sensors to efficiently detect, classify, and identify drones operating within no-drone zones. However, comprehensive research that explores various types of AAM surveillance networks, taking into account diverse requirements such as tracking different types of AAM aircraft and the use of different types of sensors to build a AAM surveillance network, has not yet received sufficient attention. While research has been conducted on state of the art sensor types related to ground based detect and avoid systems for UAS, as presented in \cite{scheff2021state}, these sensor types have not been thoroughly examined to develop different types of AAM surveillance networks. In this study, we investigate six potential sensor types suitable for AAM surveillance and assess different types of surveillance sensor networks: homogeneous sensor networks consisting of a single sensor type and heterogeneous sensor networks consisting of various sensor types. Additionally, we conduct a cost-benefit analysis of the AAM surveillance network to provide justification for investing in AAM surveillance infrastructure.

% In another study } \cite{huang2023remote}, \hl{the authors conducted a study that focused on the implementation of remote ID systems for AAM. Their objective was to assess the coverage of broadcast-receive remote ID technologies in scenarios like urban package delivery missions. Their findings revealed the approximate number of remote ID receivers required, ranging from 10 to 5000 receivers, depending on the technology and coverage goals. The authors also mentioned that their research, along with other studies on remote ID bandwidth and deployment strategies as presented in} \cite{rd1, rd2}, \hl{can provide valuable guidance to municipal authorities and AAM stakeholders for the future deployment and maintenance of remote ID systems. In} \cite{ippolito2023structurally} \hl{ }

% \cite{ippolito2023structurally}, \hl{advancements in wireless communication were harnessed to address challenges in NASA's AAM initiative, driven by the push for AAM operations in densely populated urban areas. Drawing inspiration from distributed sensing and smart spaces, where the environment integrates sensing, processing, and communication, they proposed a dynamic, topologically adaptable surveillance framework.} 

% , contributing to NASA's research in this field and preparing for eventual flight testing
% This GBVT employs image subtraction and blob detection, providing azimuth and elevation angles of detected objects from uncrewed aerial system tests. 
% \hl{In this literature review, we mostly focus on various techniques to solve optimal location selection problem. }

\subsection{Location Selection Problems and Surveillance Network Design}

The motivating application of SAND model is AAM surveillance network design which involves solving a sensor location selection problem. 
% determining the optimal locations to place AAM sensors by
% to maximize the coverage and detection performance while minimizing the cost of the network. The goal of this design is to develop an efficient AAM surveillance network that fulfills the desired coverage and detection requirements. 
In general, the location selection problem is concerned with determining the best locations for new facilities and services with respect to performance metrics such as cost, revenue, profit, travel time, distance, customer satisfaction, etc. Such problems arise in various fields such as remote sensing, geography, economics, and operations research \cite{blair1993location}. Researchers and practitioners have developed various techniques to address the location selection problem, which primarily include multi-criteria decision-making (MCDM), machine learning (ML), heuristics, metaheuristics, and mathematical optimization.

\textit{Multi-Criteria Decision-Making}: MCDM is a well-known approach used for tackling location selection problems. As the name suggests, MCDM determines the optimal location for various types of facilities --- such as a new manufacturing plant, a retail store, a hospital, a distribution center, a transportation hub, or a renewable energy facility --- based on multiple criteria or objectives. In the context of location selection problems, the criteria can include factors such as proximity to suppliers or customers, transportation costs, availability of labor, and many other factors that can affect the desirability of a location. The use of MCDM in location selection problems has been extensively studied in the literature. Three such representative papers are discussed next. The locations of manufacturing facilities were determined using MCDM in \cite{govindan2016effect}, taking into account criteria such as access to raw materials, labor force, and transportation infrastructure. This study considered several criteria such as economic, environmental, and societal factors in their facility location selection problem for sustainable development in manufacturing firms. An analytical hierarchy process (AHP) was also used here to evaluate the weights of these criteria, and a technique for order preference by similarity to ideal solution (TOPSIS) was used to rank the alternative potential locations. MCDM was also used to determine the best locations for hospitals and clinics, taking into account factors such as patient population, access to public transportation, and proximity to other healthcare providers. A location selection problem was solved in \cite{alosta2021resolving} using a MCDM approach to find the emergency medical service centers. In this study, AHP was used to determine the weights of criteria including response time, demand, coverage area, and ambulance workload. Then the different alternative locations of service centers were ranked using a technique known as ranking of alternatives through functional mapping of criterion sub-intervals into a single interval (RAFSI). A MCDM approach based on a fuzzy approach was presented in \cite{pinar2019healthcare} for determining the location of healthcare facilities. While the fuzzy logic approach can model complex systems using linguistic variables and fuzzy sets it does not however guarantee an optimal solution. MCDM techniques, such as AHP, TOPSIS, and RAFSI, rely on human judgement to evaluate and weight the criteria. There may be differences in opinion or interpretation among the decision-makers, which can lead to inconsistent or biased results, particularly if the criteria are not well-defined or if the weighting scheme is not properly calibrated. Hence, MCDM may not consistently yield the optimal locations. In contrast, mathematical optimization can consistently determine the optimal solution for the location selection problem by utilizing a mathematical model which adequately captures all relevant preferences and constraints.

% While these methods are useful in ranking and evaluating alternatives, they may not capture the full complexity of the problem and may not provide an optimal solution. 
% The fuzzy logic approach, on the other hand, uses linguistic variables and fuzzy sets to model complex systems, and the drawback of using a fuzzy approach is that it does not guarantee an optimal solution and can lead to a lack of precision in the decision-making process.

% \subsubsection{Machine Learning}
% retail store location selection problem; healthcare waste management facility selection problem; hub location problem in logistics and transportation industries 
% data mining framework, K-means algorithms with the stratified best-worst method
% benchmark: traditional heuristics

\textit{Machine Learning}: In recent years, ML techniques have been increasingly applied to solve location selection problems. Next, three papers are presented that applied ML techniques to solve logistic hub and sensor location selection problems. A ML-based algorithm framework was proposed to solve hub location problems in logistics and transportation industries. This framework consisted of a deep-learning probabilistic hub-ranker that ranks the priority of nodes to be chosen as hubs \cite{li2023machine}. To evaluate the effectiveness of this approach, the study created 11,000 small networks, each with 25 nodes, using a proposed data augmentation technique. These synthetic networks were divided into three sets: 10,000 for training models, 500 for validation during training to prevent over fitting, and 500 for model evaluation. \cite{semaan2017optimal} presented a ML-based method for optimal sensor placement in the flow over an airfoil equipped with a Coanda actuator. The method utilized a random forest algorithm to construct ML models that predicted a response function based on input data from 96 sensors measuring pressure and skin friction coefficient. The optimal sensor positions were determined by identifying the most important input variables in the ML model. A limitation of the method was its reliance on many sensors during the training phase, which made it challenging to implement experimentally. A ML algorithm was proposed and implemented in \cite{uyeh2022online} to select optimal sensor locations in controlled environment agriculture, where the macro-climate affects the micro-climate, making it challenging to predict the ideal conditions. The algorithm used temperature and humidity data from 56 different locations which were collected over a year, processed to remove outliers, and transformed to other air properties. The results showed that 3 to 5 sensors were needed, and there were similar sensor locations for different air properties. While ML has shown great promise in solving complex optimization problems such as the location selection problem, it has some drawbacks compared to mathematical optimization methods. One potential drawback is that ML models are typically designed for specific problem settings and may not be easily adaptable to other problem settings or variations, which can be a disadvantage when dealing with a new, evolving field like AAM. Though retraining ML models is possible to adapt to new problem settings, this process can be time-consuming and computationally expensive, which may limit the practical usefulness of these approaches in dynamic and rapidly evolving fields. Another challenge is that ML models require large amounts of data to train effectively, and the scarcity of AAM-related training data can hinder their applicability. In contrast, mathematical optimization methods are highly flexible and can be customized to suit a wide range of problem settings by adjusting the objective function, decision variables, and constraints to model the new problem. Moreover, unlike ML models, mathematical optimization methods provide guarantees for generating optimal solutions. Given the limited availability of AAM-related data and the need for optimal solutions, mathematical optimization methods appear to be a more suitable option for addressing location selection problems in emerging transportation sectors such as AAM.

\textit{Heuristics and Metaheuristics}: Heuristic and metaheuristic algorithms have also been used to address location selection problems due to their ability to efficiently handle complex and large-scale problems.  A proposed theory was described in \cite{dhillon2002sensor} for optimizing the placement and number of sensors in a sensor network using a greedy heuristic. The sensor field was represented as a grid of points, and the optimization framework addressed coverage optimization under constraints of imprecise detection and terrain properties. The article explained how obstacles in the terrain were modeled in the framework, and the sensor placement algorithm used a greedy heuristic to determine the best placement of one sensor at a time. The algorithm was iterative and terminates either when a preset upper limit on the number of sensors was reached or when sufficient coverage of the grid points was achieved. To solve a problem for locating fire stations and allocating resources to different stations based on dynamic traffic conditions, two metaheuristic algorithms --- particle swarm optimization and artificial bee colony --- were used in \cite{hajipour2022dynamic}. In another study, \cite{indu2009optimal} used genetic algorithm to optimize the placement of security cameras, providing maximum coverage of user-defined priority areas and minimizing the probability of occlusion due to moving objects by covering each priority area with multiple cameras. \cite{darabseh2020ads} proposed a solution for determining the optimal placement of ADS-B receivers on the ground. In \cite{darabseh2020ads}, a genetic algorithm was utilized to determine the optimal placement of ADS-B receivers on the ground in the vicinity of Frankfurt airport in Germany. The authors first identified the required number of sensors to ensure adequate coverage of the small geographical area. The algorithm was designed to search for the best local minimums, or near-optimal solutions. Then, the authors quantified the deviation of the sensor configuration generated by the algorithm from the optimal solution. \cite{inoue2023data} also implemented a metaheuristic method for solving optimal sensor placement problems using an annealing machine. One major drawback of such heuristic and metaheuristic algorithms compared to mathematical optimization is that they do not guarantee finding the global optimal solution. They involve stochastic search processes, which means that they may converge to a sub-optimal solution depending on the starting point and the algorithm parameters.

\textit{Mathematical Optimization}: Mathematical optimization is a widely used approach to solve location selection problems that aims to identify the best solution by minimizing or maximizing an objective function subject to a set of constraints. The advantage of using mathematical optimization is that it allows us to find the optimal solution that satisfies all the constraints with high precision and efficiency. \cite{chen2023developing} presented a decision-making process to select the location for public truck parking lots in urban areas using mixed-integer programming. The process included candidate location selection by spatial analysis and optimal location determination using the competitive p-median algorithm. A constrained multi-objective optimization problem with mixed-integer programming was developed in \cite{yu2022synchronous} to simultaneously determine the placement of wireless sensors and sinks that minimize energy consumption and maximize information effectiveness for structural health monitoring (SHM). \cite{ostachowicz2019optimization} provided an overview of the state of the art in the area of optimization of sensor placement for SHM applications. The optimal sensor placement problem addressed in \cite{vecherin2008optimal} aimed to select appropriate types and locations of sensors that could cover high-value terrain areas while minimizing a cost function. The probability of detection was assumed to depend on terrain conditions and obstructions. Two strategies were used for optimal sensor placement: the initial strategy utilized a heuristic and fast approach that involved placing sensors one-by-one in the location where they were most needed, while the second strategy was a binary linear programming solution that determined the global optimum of the total cost of sensors, without allowing for the sequential placement of sensors. 

% The aim was to reduce the deployment cost of SHM systems without compromising on the quality of the monitoring approach. The article coverd the optimization problem definition and each step of the optimization process, as well as the techniques and algorithms used in the literature. 

%highlights of our contribution (for location selection problem)

%revised
The location selection problem is addressed in various fields, but to our knowledge, no paper in the literature addresses the AAM surveillance sensor placement problem. After evaluating the strengths and limitations of various methods discussed above, mathematical optimization is determined to be the most suitable method for addressing the surveillance network design problem of AAM. Therefore, we use an considering the current limitations of the methods, an optimization approach has been identified as the most appropriate solution to address the gap. Hence, we develop the SAND optimization model, a binary integer linear programming (BILP) model, to design an optimized AAM surveillance network with minimum cost. In developing the SAND model, we consider the distinctive equations and parameters related to AAM, including different AAM sensor types with varying radii and other specifications, as well as the minimum required detection probability of sensors, different terrain types of a given area, and the probability of sensor detection based on terrain types. These parameters can be easily adjusted in the model's parameter setting to update the AAM surveillance network as the field evolves. Although a few articles, as mentioned above, have incorporated terrain types and sensor detection probabilities into their generic sensor location optimization models, our work stands out by utilizing real-world data and analyzing a larger number of candidate sensor locations. To achieve this, we divide the AAM surveillance area into a larger number of smaller blocks, accounting for higher computational complexity. Our approach demonstrates the ability to handle large-scale data, setting it apart from other studies. Additionally, our model can serve as a baseline for future research to compare the accuracy of solutions obtained through other methods such as ML, heuristic, MCDM, etc. Furthermore, we conduct a thorough cost-benefit analysis of AAM surveillance network and LASIC, taking into account the AAM surveillance data types and the size generated by the sensors while detecting and tracking AAM vehicles, among other essential factors that are vital to build LASIC, as detailed in Section \ref{sec3}.

% We ensure to analyze the terrain of the area and capture its actual size.

\section{Methodology}\label{sec3}

% An outline of our methodology related to AAM surveillance network and LASIC is presented in Figure \ref{fig:flow}. In this section, we first presented the SAND model and then conducted cost-benefit analysis to justify investment in SAND and LASIC infrastructure.

\begin{figure}[H]
  
\begin{adjustwidth}{-\extralength}{0cm}
%\centering %% If there is a figure in wide page, please release command \centering
  \centering
    \includegraphics[width=18cm,height=10cm]{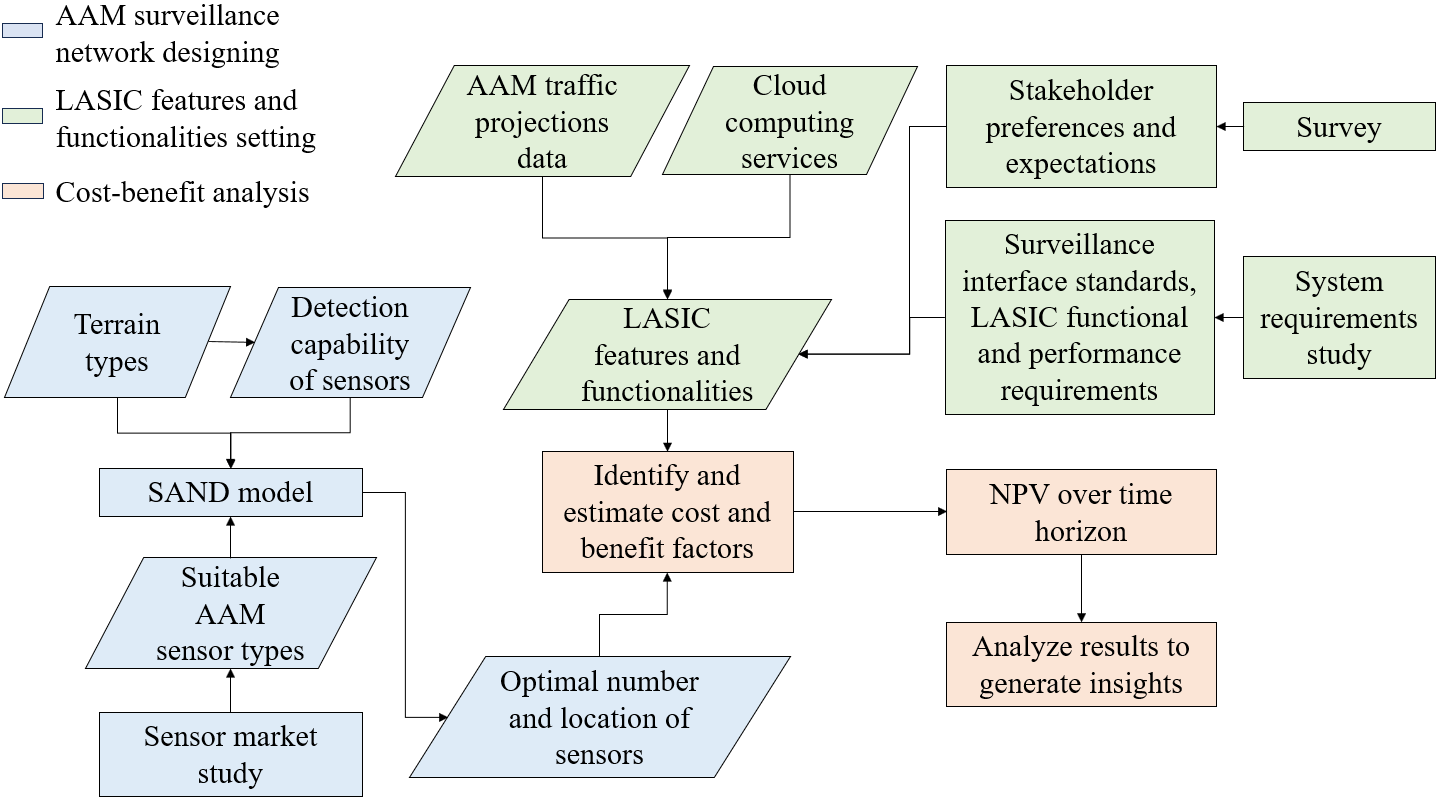}
\end{adjustwidth}
    \caption{A flow chart illustrating the steps associated with AAM surveillance network designing  and cost-benefit analysis of AAM surveillance network and LASIC.}
    \label{fig:flow}
\end{figure}

This section aims to develop the SAND model for designing an AAM surveillance network and also considers a framework for a low altitude surveillance information clearinghouse (LASIC) to receive surveillance data from the optimized AAM surveillance network. We then conduct a cost-benefit analysis to justify the investment in the AAM surveillance network-LASIC infrastructure. An outline of our methodology is presented in Figure \ref{fig:flow}, which is divided into three parts: the first part involves designing the AAM surveillance network, the second part involves setting LASIC features and functionalities, and the third part involves conducting a cost-benefit analysis of AAM surveillance network and LASIC. The blue boxes represent the inputs of the SAND model, including AAM sensor types, terrain types, and sensor detection probabilities, as well as its outputs, such as the optimal number and location of sensors, which are further discussed in Section \ref{SND}. The green boxes correspond to the inputs LASIC feature and functionalities, including AAM traffic projections data, cloud computing services, and insights from survey and system requirements study that are elaborated upon in Sections \ref{ccc} and \ref{cf2}. The orange boxes represent the steps associated with the cost-benefit analysis of the AAM surveillance network and LASIC. The results of the first and second parts are utilized to identify and estimate the cost and benefit factors, which are discussed in detail in Section \ref{CBA} and Section \ref{benefit_fact}, respectively. The output of the cost-benefit analysis is the NPV over the analysis period, which is explained in Section \ref{NPV}.

% To undertake a 10-year city-based assessment of LASIC for the State of Ohio, we computed net present value (NPV) of AAM surveillance network and LASIC for the six major cities of Ohio (SMCO): Columbus, Cleveland, Cincinnati, Akron, Toledo, and Dayton. The NPV is a measure of the future return on investment expected from an investment in a project in terms of today's dollars. The NPV metric takes into account the time value of money and future cash flows, which is further discussed in Section \ref{NPV}. The formation of SMCO was predicated on the finding of significant demand potential for AAM use cases in those cities considering socioeconomic factors, such as population, population density, gross domestic product, median per capita income, cost of living, city total area, cities in motion index, human capital, etc \cite{del2021infrastructure}.

% We carried out the cost-benefit analysis into two major steps: 1) identification and evaluation of cost and benefit factors and 2) calculation of net present value. 

\subsection{Surveillance Sensor Network Design for Advanced Air Mobility} \label{SND}

We construct a BILP model to solve the optimal sensor placement problem for designing the surveillance network for AAM. The objective of this study is to design a network of surveillance sensors to track AAM aircraft, flying in the lower altitude over a city, with a minimum sensor cost satisfying two constraints: 1) sensor(s) at a location must provide a minimum required detection probability and 2) area across a city must be fully covered by the network. Given these restrictions, the goal is to determine the optimal locations for sensor placement and the number of sensors needed to be placed within a city to detect and track the aircraft, where the objective function focuses on minimizing the total sensor cost. 

This section highlights the crucial terms that must be taken into consideration while formulating the SAND model, such as various types of AAM surveillance sensors, detection probability of a sensor and the impact of terrain types on the sensor detection probabilities. Subsequently, we present the formulation of the SAND model by modeling the relevant surveillance area where the model is implemented. A summary of the notations used in designing the SAND model is presented in Table \ref{notation}.

\begin{longtable}{@{}p{3cm}p{10.4cm}@{}}
\caption{Parameters, indices, and decision variables in the SAND model.} \label{notation} \\
\toprule
\textbf{Parameters} & \textbf{Definition} \\ \midrule
\endfirsthead
% \multicolumn{2}{c}{\tablename\ \thetable\ (continued)} \\
% \toprule
% \textbf{Indices} & \textbf{Definition} \\ \midrule
% \endhead
$M$ & A rectangular mesh. \\
$F$ & Transformation function of GCS coordinates to PCS coordinates. \\
$\lambda_p$ & Longitude of the $p$-th point in GCS. \\
$\phi_p$ & Latitude of the $p$-th point in GCS. \\
$n_a$ & Number of points along the x-axis of $M$. \\
$n_b$ & Number of points along the y-axis of $M$. \\
$L_a$ & Length of the area along the horizontal axis. \\
$L_b$ & Length of the area along the vertical axis. \\
$L$ & Block side length. \\
$\rho$ & Range of a sensor. \\
$P^M$ & Set of all points in $M$. \\
$Z^M$ & Set of all blocks in $M$. \\
$T$ & Set of terrain types associated with each block in $Z^M$. \\
$S$ & Set of potential sensor types. \\
$\omega^s_T$ & Detection probabilities for all combinations of terrain types in $T$ and sensor types in $S$. \\ 
$R_s$ & Sensor range for a sensor of type $s$ in $S$.\\
$\omega^s_z$ & Probability of detecting an AAM aircraft with a sensor of type $s$ on block $z$. \\
$I(z)$ & Indicator function that equals 1 if block $z$ belongs to the area, and 0 otherwise. \\
$Q$ & Number of blocks removed from $Z^M$. \\
$C$ & Set of center points of blocks in $Z$. \\
$d_{e,i}$ & Euclidean distance between a sensor location $e$ in $C$ and a point $i$ in $P^M$. \\
$A^s_e$ & Set of coordinates of the points covered by a sensor of type $s$ at location $e$. \\
$B^s_e$ & Set of blocks covered by a sensor of type $s$ at location $e$. \\
$\zeta^s_e$ & Mean of all the probability of sensor detection values for blocks in $B^s_e$ for a sensor of type $s$. \\
$\tau^s_e$ & Probability of misdetection of a sensor of type $s$ at location $e$. \\
$r$ & Minimum required detection probability.\\ 
$\kappa^s_e$ &  Number of independent sensors of type $s$ required to achieve a minimum required detection probability at location $e$. \\
$\psi^s$ & Cost of a sensor of type $s$. \\ 
$\delta^s$ & Number of sensors needed for sensor type $s$ to provide $360^\circ$ coverage at a location. \\
\midrule

\textbf{Indices} &  \\ \midrule
$p$ & $p$-th point in GCS. \\
$j$ & $j$-th row of blocks in $M$. \\
$k$ & $k$-th column of blocks in $M$. \\
$z$ & $z$-th block in $Z^M$ and $Z$.\\ 
$e$ & $e$-th candidate sensor location in $C$. \\
$i$ & $i$-th point in $P^M$. \\
$s$ & $s$-th sensor type in S.\\

\midrule

\textbf{Decision Variables} &  \\  \midrule

$\lambda^s_e$ & Binary variable representing whether a sensor of type $s$ is placed at location $e$. \\
$\gamma_z$ & Binary variable representing whether block $z$ is covered by at least one sensor. \\
\bottomrule
\end{longtable}

\subsubsection{Surveillance Sensor Types} \label{sensortype}

After studying the existing sensor market, six types of sensors are deemed to be suitable for AAM traffic surveillance. A brief overview of these surveillance sensor types are provided in this section. 

% We focused on the suitable AAM sensors identified in the aircraft surveillance source study.
% we selected one sensor from each sensor type in our analysis. The sensors from six sensor types are described below. 

\noindent 1. Radar: Both cooperative and non-cooperative aircraft can be detected and tracked using ground-based radars. The radar transmits electromagnetic waves signal towards aircraft which bounce off the aircraft and create a detailed image of its size, shape, and location. The radar cross-section (RCS) signature of each aircraft type is distinctive, which leads to varying reflection patterns of radio waves. The radar utilizes these patterns to identify the aircraft type and determine its position, velocity, and travel direction. 

\noindent 2. Automatic Dependent Surveillance-Broadcast: Automatic Dependent Surveillance-Broadcast (ADS-B) is a surveillance system that allows an aircraft to periodically broadcast and track its location via satellite navigation. Currently, FAA  acknowledges ADS-B as a key enabler for trajectory-based air traffic management in the future. 

\noindent 3. Remote ID: The ability of sUAS and eVTOL to broadcast identification and location data during its flight is known as \textit{remote identification (remote ID)}. 

\noindent 4. Radio Frequency: Like the radar, RF sensor is also able to accurately detect and categorize aircraft. However, RF sensors can detect and track small drones that may not be detectable by radar, particularly at low altitudes where the radar signal may not reflect off the drone as effectively as it would off a larger aircraft. Also, RF sensors can be more effective than radar in urban or cluttered environments where there may be many buildings, trees, and other obstacles that can reflect or absorb radar signals. RF sensors are less affected by these obstacles because their signals can penetrate walls and other structures, making them useful for monitoring drones in indoor or urban environments. The key advantages of the RF sensor system include its low cost, ease of installation, and simplicity of integration with several other sensors, including cameras and radars. 

\noindent 5. Acoustic: An audio pattern that is transmitted by an aircraft's propeller can be detected by acoustic sensors and used for aircraft positioning and classification. It uses passive acoustic sensor technology with no RF emissions, where the solid state-sensor is an array module including digital microphones and digital processors.

\noindent 6. Electro-Optical/Infrared Camera: An Electro-Optical/Infrared (EO/IR) system is a type of electronic equipment that combines electro-optical and infrared sensors to offer accurate optical information of air traffic in the airspace within its coverage range at any time. EO/IR systems can be used to carry out object tracking, assess threats from a certain distance, or monitor other aircraft or ground obstructions that must be avoided. 

\subsubsection{Detection Probability of Sensors and Terrain Types} \label{detect_prob}

Detection probability is a crucial performance parameter for sensors, representing the probability of a sensor detecting an object within its field of view. Terrain types of an area, such as hill, open space, and water bodies, can significantly affect the detection probability of a sensor by obstructing its field of view and reducing its detection range. Different terrain types can have a range of effects on the sensor detection probabilities, with some terrain types having a greater impact than others. For example, hilly terrain types can have a substantial impact on the detection probability of sensors due to their obstruction, while open spaces can provide ideal conditions for achieving higher sensor detection probability values. To obtain the sensor detection probability values for sensors, there are two primary methods. The first method involves using manufacturer-provided data based on laboratory testing and simulations. Sensor manufacturers can provide data of detection probability for different terrain types in sensor data sheets. This data can be used to estimate the sensor's performance in various environments. The second method involves testing sensors on-site in the specific terrain type.

% for Advanced Air Mobility

% The objective of this study is to design a network of surveillance sensors to track AAM aircraft, flying in the lower altitude over a city, with a minimum sensor cost satisfying two constraints: 1) sensor(s) at a location must provide a minimum required detection probability and 2) area across a city must be fully covered by the network. Given these restrictions, the goal is to determine the optimal locations for sensor placement and the number of sensors needed to be placed within a city to detect and track the aircraft, where the objective function focuses on minimizing the total sensor cost. 

% In the context of this problem, 
% in which the surveillance network needs to be set up

% \subsubsection{Modeling the Surveillance Area} \label{model_area}

\subsubsection{SAND Model Formulation} \label{model}

In developing the SAND model, the process of modeling a surveillance area is critical in designing an effective AAM surveillance network, as it directly impacts the placement of sensors. This involves identifying the surveillance area and determining feasible sensor locations, while also considering that not all locations within the surveillance area may be feasible, such as those located over water bodies. Through this process, we can optimize the placement of sensors by accurately capturing the terrain types of an area and assessing their impact on sensor placement and detection probabilities. 

To start modelling the surveillance area, first, a rectangular mesh denoted by $M$ is used to divide a given area which needs AAM surveillance into a set of points and small square blocks. To create such a mesh, the city is first overlaid by a rectangle defined by four geographic coordinate system (GCS) points on a world map, labeled $A$, $B$, $C$, and $D$. However, the use of latitude and longitude to define locations in GCS coordinates means that distances on the Earth's surface can result in blocks of unequal size and non-parallel mesh lines in $M$. To address this issue, it is necessary to use a projected coordinate system (PCS) that maps the Earth's surface to a 2D Cartesian plane when creating on a mesh for an area. This ensures that $M$ is created with equal-sized blocks and straight, parallel mesh lines, which is important in the optimal sensor placement problem for accurate measurements and simplified visualization and analysis. In Equation \ref{py}, $F$ represents the transformation function of GCS coordinates to PCS coordinates, and $\lambda_p$ and $\phi_p$ represent the longitude and latitude, respectively, of the $p$-th point in GCS. The output of the transformation is the corresponding $(x_p, y_p)$ PCS coordinates.
 
% To convert the four GCS points $A$, $B$, $C$, and $D$ into PCS coordinates, the 'pyproj' package in Python is used
  
\begin{equation} \label{py}
% (x_p, y_p) = pyproj.transform(proj_{GCS}, proj_{PCS}, \lambda_p, \phi_p), \indent p \in \{A, B, C, D\}.  
(x_p, y_p) = F(\lambda_p, \phi_p), \indent p \in \{A, B, C, D\} 
\end{equation} 

% $proj_{GCS}$ represents the projection method used for the GCS, while $proj_{PCS}$ represents the projection method used for PCS.

% A 2D mesh, denoted by $M$, is generated from $A$, $B$, $C$, and $D$ using parameters $n_a$ and $n_b$, where $n_a$ and $n_b$ respectively specify the number of points along the x-axis and y-axis. The values of $n_a$ and $n_b$ in Equation \ref{n,S,R} are determined by the length of the area along the horizontal axis, $L_a$, and the length of the area along the vertical axis, $L_b$, respectively, as well as by the block length $L$:

% The parameters $n_a$ and $n_b$ in Equation \ref{n,S,R} specify the number of points along the x-axis and y-axis, which are determined by the length of the area along the horizontal axis, $L_a$, and the length of the area along the vertical axis, $L_b$, respectively, as well as by the block side length $L$.

% \begin{equation}\label{n,S,R}
% n_a = \left\lceil\frac{L_a}{L}\right\rceil, \quad n_b = \left\lceil\frac{L_b}{L}\right\rceil, \quad \textrm{subject to } R \geq \frac{L}{\sqrt2}
% \end{equation}

\begin{equation}
n_a = \left\lceil\frac{L_a}{L}\right\rceil, n_b = \left\lceil\frac{L_b}{L}\right\rceil, \textrm{subject to } \rho \geq \frac{L}{\sqrt2},
\end{equation}

\noindent where $n_a$ and $n_b$ specify the number of points along the x-axis and y-axis. These parameters are determined by the length of the area along the horizontal axis, $L_a$, and the length of the area along the vertical axis, $L_b$, respectively, as well as by the block side length $L$.
% \noindent If $L$ is larger than the sensor coverage range denoted by $R$, there may be areas within the block that are not covered by any sensors, leading to blind spots or gaps in the coverage that can compromise the effectiveness of the AAM surveillance network. Therefore, the condition that $L<\rho$ ensures that each block is adequately covered by the sensor(s). 
The diagonal of a square block is represented by $\sqrt2L$, with the center point of a block being considered as the potential location for placing a sensor. Hence, the distance from the sensor location to a corner point of the block is $L/\sqrt{2}$. If this distance is greater than the sensor coverage range denoted by $\rho$, the block will not be covered by the sensor. Thus, $\rho \geq \frac{L}{\sqrt2}$ ensures that each block is adequately covered by the sensor(s) and to avoid blind spots or gaps in the coverage of the AAM surveillance network. The resulting mesh $M$ consists of $n_a \times n_b$ points, denoted by $P^M$, where $(x_i, y_i)$ corresponds to the $i$-th point on the 2D map. Therefore, we can express the set of all points in $M$ as

\begin{equation}
P^M = \{(x_1, y_1), (x_2, y_2), \ldots, (x_{n_a \times n_b}, y_{n_a \times n_b})\}.
\end{equation}

\noindent The set of blocks within $M$ is defined by

\begin{equation}\label{ZM}
Z^M = \left[\left[P^M_{j \times n_a + k}, P^M_{j \times n_a + k + 1}, P^M_{(j+1) \times n_a + k}, P^M_{(j+1) \times n_a + k + 1}\right]\right]_{j=1, k=1}^{n_b-1, n_a-1},
\end{equation}

\noindent where the number of blocks in $M$ is $(n_a-1) \times (n_b-1)$. For each adjacent point pair in $M$, a block with four corner points is created using  $j$ and $k$, the indices of $n_a$ and $n_b$, respectively, such that $1 <= j < (n_b-1)$ and $1 <= k < (n_a-1)$.

The probability of detecting an AAM aircraft flying over a given block by a given sensor type is determined by the terrain type of that block. Let $T$ be a set of terrain types associated with each block within the mesh for a given area, and $S$ be a set of potential sensor types. The detection probabilities for all combinations of terrain types and sensor types is represented by the matrix $\omega^s_T$. Let $z$ be an index used to iterate through $Z^M$, where $z$ ranges from 1 to $(n_a-1) \times (n_b-1)$. $T_z$ represents the terrain type of the $z$-th block in $T$, and $s$ is the index for the sensor type set $S$. The probability of detecting an AAM aircraft with a sensor of type $s$ on block $z \in Z^M$, where block $z$ has terrain of type $T_z$, is given by

\begin{equation}\label{omega}
\omega^s_z = \omega^s_{T_z},  \indent \forall z \in Z^M, \forall s \in S. 
\end{equation}

% The SAND model
% To ensure that sensor locations are within the area boundaries by excluding the outer blocks of the mesh that do not belong to the area. This is achieved by setting the detection probability of sensors in these outer blocks to zero and removing them from the block list $Z^M$.

As the surveillance area of interest will likely have an irregular shape, some of the outer blocks of the rectangular mesh will not belong to the area. These outer blocks are removed from the block list $Z^M$.  The remaining set of blocks present within the area is represented as $Z = {z \in Z^M \mid I(z) = 1}$, where $I(z)$ is an indicator function that equals 1 if the block $z$ belongs to the area, and 0 otherwise. %The detection probability of sensors in the blocks outside the area is set to zero, i.e., $\omega^S(b) = 0$ for $b \in Z^M \setminus Z$ and for a sensor of type $s$. 
The number of blocks that are removed from $Z^M$ is given by 
\begin{equation} \label{Q}
Q = \left|\left\{ z \in \mathbb{Z}^M \mid I(z) = 0 \right\}\right|.
\end{equation}

\noindent Then, $Z$ is the updated block list, where $z$ ranges from 1 to $[{(n_a-1) \times (n_b-1)}-Q]$. By doing so, the model can approximate the actual shape of the area and select optimal sensor locations within the area. The set of center points of blocks in $Z$ is considered as the candidate sensor location. The coordinates of candidate sensor location is represented by 

\begin{equation} 
C = \{(\alpha_e, \beta_e) \mid e \in Z\}, 
\end{equation}

\noindent where a potential sensor location is denoted by $e$. Since water blocks cannot be selected as sensor locations, the center points of water blocks are excluded from $C$ to ensure they are not considered as potential sensor locations. Then, the Euclidean distance, 

\begin{equation}
d_{e,i} = \sqrt{(\alpha_e - x_i)^2 + (\beta_e - y_i)^2}, \indent \forall e\in C, \forall i \in P^M, 
\end{equation}

\noindent is computed to measure distance between a sensor location $e \in C$ and a point $i \in P^M$. If a sensor of type $s$ is placed at $(\alpha_e, \beta_e)$, it can cover the points that are within its sensor range $R_s$. To find the points within the range $R_s$ of a sensor of type $S$ placed at $(\alpha_e, \beta_e)$, the set of coordinates of the points  covered by the sensor is given by

\begin{equation}
A^s_e= \{i \in P^M \mid d_{e,i} \leq R_s\}, \indent \forall  e \in C,  \forall s \in S. 
\end{equation}

To determine the set of blocks that would be covered by each sensor at each candidate location, let $B^s_e$ be the set of blocks that would be covered when sensor $s \in S$ is placed at location $e$. Algorithm \ref{algo1} iteratively checks whether each point $o$ in the block $z$ is in the set $A^s_e$. If any point in $z$ is not in $A^s_e$, the algorithm sets the value of $A^s_e$ to False for that block $z$. If all the points in $z$ are in $A^s_e$, the algorithm indicates that a sensor of type $s$ at location $e$ covers the block $z$.

% \begin{equation} \label{B^s}
% \forall z \in Z, \forall o \in Z_z, (o \in A^s_e) \rightarrow B^s_e\leftarrow B^s_e \cup {z},  \indent \forall e \in C \forall s \in S
% \end{equation}

\begin{algorithm}[H]
\caption{\textls[-10]{Computing the set of blocks covered by each sensor at each candidate location.}} \label{algo1}
%\SetAlgoLined
% \KwResult{Set of blocks that would be covered when each sensor is placed at each candidate location}
\textbf{for} each sensor of type $s$ in $S$ \textbf{do}

\hspace{2em}\textbf{for} each location $e$ in $C$ \textbf{do}
initialize an empty set $B^s_e$;

\hspace{4em}\textbf{for} each block $z$ in $Z$ \textbf{do}
all points in $z$ are in $A^s_e$ = True;

\hspace{6em}\textbf{for} each point $o$ in block $z$ in $Z$ \textbf{do}

\hspace{8em}\textbf{if} $o$ is not in the set $A^s_e$ \textbf{then}
all points in $z$ are in $A^s_e$ = False;
break;

\hspace{8em}\textbf{if} all points in $z$ are in $A^s_e$  \textbf{then}
add $z$ to the set $B^s_e$.
\end{algorithm}

% \begin{equation} \label{B}
% B^s_e = \big\{z \in Z \mid \forall\, q_{i,j} \in Z_{z_i},\, q_{i,j}\in A_e\big\}, \indent \forall e \in C \forall s \in S
% \end{equation}

% In Equation \ref{pd}, $\omega^s_z$ represents the probability of sensor detection for block $z$ for a sensor of type $s$, and $\zeta^s_e$ which is a mean of all the probability of sensor detection values for blocks $z$ in $B^s_e$ for that sensor type. 
% The notation $\frac{1}{|B^s_e|}$ represents the fraction of the total number of blocks in $B^s_e$, which serves as a normalization factor to ensure that the final probability value is within the range of 0 to 1. 

To consider the average detection probability of a sensor placed at $e \in C$, a mean of all the probability of sensor detection values for blocks $z$ in $B^s_e$ is computed for a sensor of type $s$, which is denoted by

\begin{equation}\label{pd}
\zeta^s_e = \frac{1}{|B^s_e|}\sum_{z \in B^s_e}\omega^s_z, \indent \forall e \in C, \forall s \in S.
\end{equation}

The probabilistic framework of sensor detection probability presents an important consideration to improve the detection of aircraft in the airspace by understanding the probability of detection. When an aircraft is present in the airspace, the likelihood that a sensor will detect it is known as the probability of detection. On the other hand, the probability of misdetection refers to the likelihood of not detecting an aircraft when it is actually present. An effective solution to address this issue is to ensure sufficient sensor coverage at a specific location by installing an adequate number of sensors. This ensures that at least one sensor can track aircraft that meet the minimum detection probability requirement, which is the minimum probability of detecting an aircraft that must be met by a sensor to ensure reliable detection. %To solve this issue, the installation of enough sensors is necessary to ensure sufficient coverage at a specific location. This way, at least one sensor can track aircraft that meet a minimum detection probability requirement.% % \cite{vecherin2008optimal}
For example, if two sensors are placed at a location, each with a detection probability of 0.8, the probability that at least one of the sensors will detect the aircraft is 0.96. Hence, it is essential to consider the probability of misdetection of a sensor at location $e$ for a sensor of type $s$, which is denoted by

\begin{equation}\label{tau}
\tau^s_e = 1 - \zeta^s_e, \indent \forall e \in C, \forall s \in S.
\end{equation}

To achieve a minimum detection probability requirement, denoted by $r$, the number of independent sensors required is needed to determine by

\begin{equation}\label{kappa}
\kappa^s_e = \frac{\log(1 - r)}{\log(\tau^s_e)}*\delta^s, \indent \forall e \in C, \forall s \in S
\end{equation}

\noindent for a sensor of type $s$ at location $e$. $\delta^s$ represents the number of sensors needed for sensor type $s$ to provide $360^\circ$ coverage at a location, which depends on the field of view of that sensor type. For example, if the sensors have a field of view of $90^\circ$, then four sensors of the same types are needed to be positioned at equal intervals around the location to provide full coverage. If the sensors have a wider or narrower field of view, fewer or more sensors may be needed to ensure complete coverage, respectively.

% For each sensor of type $s$ in $S$, the corresponding $B^s$ and $\kappa^s$ values are aggregated into $B$ and $\kappa$, respectively, which are then considered as input variables for the SAND model as follows: 
% $\displaystyle B = \bigcup_{s\in S} B^s$ and $\displaystyle \kappa = \bigcup_{s\in S} \kappa^s$.

% \subsubsection{SAND Model Formulation} 

%The objective of the SAND model is to minimize the function $\theta$, as shown in Equation \ref{theta}.
The objective of the SAND model, which follows the modeling approach described above, is to minimize the function 

\begin{equation}\label{theta}
min(\theta) = \sum_{\forall s \in S}\sum_{e \in B^s}\lambda^s_e*\kappa^s_e*\psi^s
\end{equation}

\noindent by determining the optimal locations for sensor placement and the number of sensors required to achieve a desired level of detection probability, while minimizing the total sensor cost. This function depends on two parameters:  $\psi^s$ which denotes the cost per unit of a sensor of type $s$, $\kappa^s_e$, and a binary decision variable

\begin{equation}\label{lam}
\lambda^s_e =  \left\{
        \begin{array}{ll}
            1,  \quad &\text{if sensor is placed at}\, e \in C, \text{where}\, s \in S \\
            0, \quad &\text{otherwise.} \,
        \end{array}
    \right.
\end{equation}

\noindent We introduce another binary decision variable

\begin{equation}\label{gam}
\gamma_z =  \left\{
        \begin{array}{ll}
            1,  \quad &\text{if block}\, z \in Z\, \text{is covered by at least one sensor} \\
            0, \quad &\text{otherwise.} \,
        \end{array}
    \right.
\end{equation}

\noindent Constraint \ref{cons} ensures that each block $z \in Z$ is covered by at least one sensor $e$ in $C$. 

\begin{equation}\label{cons}
\sum_{e \in B^s, z \in B^s_e} \lambda^s_e \geq \gamma_z,  \indent \forall z \in Z, \forall s \in S
\end{equation}

% Several assumptions were incorporated into our study:
% tall buildings
% Secondly, we assume that all AAM flights take place within an altitude of 400 meters, with consistent sensor performance irrespective of the precise aircraft paths.
% we consider a uniform distribution of AAM flight hours across the analyzed region

Several assumptions were incorporated into our model. Firstly, we assume that there is no potential sensor obstructions posed by natural or human-made structures in the surveillance region. Secondly, we assume that all AAM flights take place within an altitude of 400 meters, which is within the range of all the different sensor types considered. Thirdly, sensor detection performance does not depend on the precise aircraft paths. Fourthly, we assume that there is no effect of weather on sensor performance. We do not consider the effect of sensor failures on the AAM surveillance network. Lastly, we assume a uniform distribution of AAM flight operations across the analyzed region.

% A number of assumptions have been made during our study:

% \begin{enumerate}

%   \item No obstruction of sensors by tall buildings has been considered.
%   \item We assumed that all AAM flights occur within a range of 400 meters or 2000 feet.
%   \item The study assumes that sensor performance remains consistent regardless of the specific aircraft trajectories.
%   % \item An assumption has been made that the growth in the number of subscribers follows the growth rate of the AAM market.
  
% \end{enumerate}

\subsubsection{Solution Algorithm}

Following the methodology described in Section \ref{SND}, the SAND model is implemented in Python 3 using the Gurobi-Python API and the model is solved using the Gurobi 10.0.1x64 Linux on a computer with 3.00 GHz × 36 Intel® Core™ i9-10980XE processor and 128 GiB of memory. The Gurobi optimizer performs a branch and bound search to find a global solution. This is a systematic technique for solving optimization problems by recursively partitioning the space into smaller branches, and then solving each branch independently. To convert the GCS coordinates into PCS coordinates for modeling the surveillance area, we use the 'pyproj' Python package which can conduct geodetic calculations and cartographic transformations \cite{pyproj}.

% \makeatletter
% \renewcommand\paragraph{\@startsection{paragraph}{4}{\z@}%
% % display heading, like subsubsection
%                                      {-3.25ex\@plus -1ex \@minus -.2ex}%
%                                      {1.5ex \@plus .2ex}%
%                                      {\normalfont\normalsize\bfseries}}
%  \setcounter{secnumdepth}{4}
% \makeatother

\subsection{Low Altitude Surveillance Information Clearinghouse Features and Functionalities} \label{ccc}

% To determine the stakeholder preferences and expectations on LASIC features and functionalities, a survey of AAM stakeholders was carried out. Based on the survey, expected LASIC features include access to real-time coverage maps, live surveillance feeds and offline archival surveillance data, and support for querying and analyzing surveillance data. The relevant functional and performance requirements of LASIC ---- namely, surveillance interface standards, sensor data sizes, and ping rate --- are determined based on a system requirements study of LASIC.

% In order to facilitate the successful functioning of LASIC, the surveillance data generated by sensors must be safely processed and stored, taking into account the features identified in the survey and system requirements study. This data can be processed either in a cloud-based server or locally-owned servers. 

% features obtained from the survey findings and system requirements study. 

The surveillance data generated by the optimized surveillance sensor network must be safely processed and stored, taking into account the stakeholder preferences and expectations on LASIC features and functionalities. To determine these preferences and expectations, a survey of AAM stakeholders was carried out. Based on the survey, expected LASIC features include access to real-time coverage maps, live surveillance feeds and offline archival surveillance data, and support for querying and analyzing surveillance data. The relevant functional and performance requirements of LASIC ---- namely, surveillance interface standards, sensor data sizes, and ping rate --- are determined based on a system requirements study of LASIC.

The surveillance data can be processed either on cloud or locally owned servers. Given the features of LASIC, a cloud based server is considered to be more suitable for hosting the surveillance data of LASIC due to the following reasons. Firstly, local servers require a large amount of time and effort to set up and maintain it. They also require a lot of space and expensive hardware. On the other hand, cloud computing servers can be a cost-effective solution for LASIC, as they eliminate the need to invest in expensive hardware and infrastructure. Instead, the LASIC operator (e.g., a given state's department of transportation) would pay only for the resources they use in the cloud, which can help build a more cost-effective LASIC in the long run by reducing operational expenses and avoiding the capital expenditures associated with maintaining and upgrading local servers. Secondly, cloud computing servers provide a higher level of security compared to local servers, as cloud computing servers invest heavily in security measures such as firewalls, encryption, and intrusion detection systems. The data gathered from surveillance must be protected in LASIC from unauthorized access, cyber threats, and breaches to ensure the privacy of the system. Any unauthorized access to the data can result in potential harm to the system and damage to the reputation of the LASIC program. The advanced security measures of cloud computing servers can ensure the confidentiality, integrity of AAM surveillance network and LASIC, which is essential to maintain the trust of their constituents and regulatory compliance. Lastly, cloud computing servers offer the advantage of being able to easily scale up or down based on changing demands. This is especially crucial in the context of LASIC and AAM surveillance as the demand for advanced air mobility is rapidly increasing. As more AAM vehicles take to the skies, the amount of data generated by these vehicles will also increase, and the computing resources required to process and store this data will need to be adjusted accordingly. Conversely, local servers have a fixed number of resources and require additional hardware investments to accommodate additional demands. This can be a significant disadvantage for LASIC operators who may need to invest in new hardware to accommodate increased demand, resulting in higher upfront costs.

Given the benefits discussed above, we deem a cloud-based server to be the most suitable host for LASIC. The cost of operating the server or utilizing cloud computing services depends on AAM traffic projections, surveillance data generated in a specific area, surveillance data types, interface standards, data sizes, ping rate, cloud components, and its pricing policies, which are discussed next.

\subsubsection{Surveillance Data Types and Sizes} \label{sec:data_types_sizes}

The cost of cloud computing services is associated with the types and sizes of surveillance data, as cloud vendors generally employ a billing model based on the number of data held or processed in the cloud infrastructure. As such, larger volumes of data processed or stored in the cloud result in higher cloud computing cost. The type and size of surveillance data generated by the sensors in the AAM surveillance network depend on the type of service level provided by LASIC and the surveillance data interface standard. The three possible service levels that can be provided to subscribers of LASIC are: informational only, radio location quality, and radio navigation quality. Of all the service levels, the radio navigation quality requires stringent data requirements that specify strict and precise data specifications, as it supports tactical deconfliction services. This is necessary to provide highly accurate navigation and positioning information, which is crucial for avoiding collisions and ensuring effective coordination between AAM aircraft during flight. To specify the interface standard for LASIC, the All-purpose STructured EUROCONTROL Surveillance Information eXchange (ASTERIX), as mentioned in \cite{129, 021, 062}, is used in this study, which is a collection of interface definitions and documentation outlining the data format standards used for transmitting a range of surveillance data.    

% These service levels have different data requirements.

% The size of the yearly total surveillance data that would be generated in a given area can be determined based on the projected yearly flight hours of AAM traffic in that area estimated for potential AAM use cases --- including passenger and cargo transportation, bridge inspections and medical items delivery  by sUAS --- and projected yearly flight hours of general aviation traffic \cite{GA} and size of surveillance data messages, the sensor ping rate 1 Hz, which is the minimum sensor data rate needed to provide real time surveillance  generated by the sensors. 

% rephrased
The size of the total yearly surveillance data generated in a given area are determined based on several factors, including the projected yearly flight hours of AAM traffic in the area for potential use cases such as passenger and cargo transportation, bridge inspections, small package delivery and medical item delivery. We obtain the data of yearly estimated AAM passenger and cargo traffic from \cite{del2021infrastructure} and conducted forecasting to estimate demand for other AAM use cases in \cite{dulia2021benefits}. Our assumption for the distribution of AAM flight hours over the given area is uniform. The size of the total yearly surveillance data also takes into account the size of the surveillance data messages, which are the packets of data generated by the sensors used for surveillance. These packets may contain information such as images, video, location data, and other types of sensor data. Additionally, the size of the total yearly surveillance data is calculated based on a ping rate of 1 Hz, which refers to the frequency at which the surveillance data is transmitted. A data ping rate of at least 1 Hz is necessary to provide real-time surveillance generated by the sensors. The AAM traffic is considered to comprise three main types of aircraft: cooperative manned aircraft, cooperative uncrewed aircraft, and non-cooperative aircraft. The surveillance message sizes associated with these types of aircraft are computed using their corresponding interface standards as defined in ASTERIX (\cite{129,021,062}), which are listed in Table \ref{tab:data_type}.

\begin{table}[H]
\caption{\label{tab:data_type} Types and sizes of surveillance data.}
\newcolumntype{C}{>{\centering\arraybackslash}X}
\begin{tabularx}{\textwidth}{lcccccc}
\toprule
% & Transition& & \multicolumn{2}{c}{}\\\cmidrule{2-2}
\textbf{Aircraft Type} && \textbf{Interface} & \textbf{Number of} & \textbf{Message} \\ 
                       && \textbf{Standard}  & \textbf{Data Items} & \textbf{Size (Bits)}  \\
                       \midrule
Cooperative Manned Aircraft && ASTERIX CAT-021 & 42 & 1136\\
Cooperative Uncrewed Aircraft && ASTERIX CAT-129 & 14 & 432\\
Non-Cooperative Aircraft && ASTERIX CAT-062 & 27 & 2648\\
\bottomrule
\end{tabularx}
\end{table}

\subsubsection{Cloud Components}

% \subsubsection{Cost Factor 2: Cloud Computing Cost}

% The surveillance sensors are envisioned to track the aircraft movement and generate surveillance data that will be ingested into LASIC for processing in real time. The cloud computing cost depends on the total yearly surveillance data generated in SMCO, which was estimated based on the AAM traffic demand, surveillance data types, associated interface standards and data sizes (as described in Section \ref{sec:data_types_sizes}). To evaluate the cost of the cloud components in this study, the Microsoft Azure Web cloud computing services and their pricing policies were considered based on the number of units of data published and received by LASIC. The Microsoft Azure provides a range of cloud-based services that can be utilized to create a platform for real-time analysis of live surveillance data. A real time analytics on big data architecture would need to be created to enable the data flowing through LASIC, which would consist of six components 1) Azure Event Hubs, 2) Azure Synapse Analytics, 3) Azure Data Lake Storage, 4) Azure Cosmos Database (DB), 5) Azure Analysis Services, and 6) Azure Power BI \cite{azure}. 
% The cloud computing cost depends on several factors, including the total yearly surveillance data generated in SMCO, surveillance data types and associated interface standards and data sizes, as discussed in Section \label{sec:data_types_sizes}, and the required cloud computing tools. 

The cloud computing cost is determined by the pricing policies of the cloud server chosen to host LASIC and the cloud components needed to enable the desired real-time and offline LASIC features and functionalities. The Microsoft Azure Web cloud computing services is considered in this study to estimate the cloud computing cost. Microsoft Azure provides a range of cloud-based services that can be utilized to create a platform for real-time analysis of live surveillance data. It can be used to continuously ingest and process LASIC data in near-real time and store the data for data archival, dissemination, querying and analytics. 

% As the pricing policies of ESP are not available yet, we considered Microsoft Azure Data Lake Storage's pricing policy to evaluate the data storage cost in our analysis.

% The ODOT event streaming platform (ESP) is considered to be used for data archival and dissemination.

\begin{figure}[H]
   
\begin{adjustwidth}{-\extralength}{0cm}
%\centering %% If there is a figure in wide page, please release command \centering
 \centering
   \includegraphics[width=15cm,height=6cm]{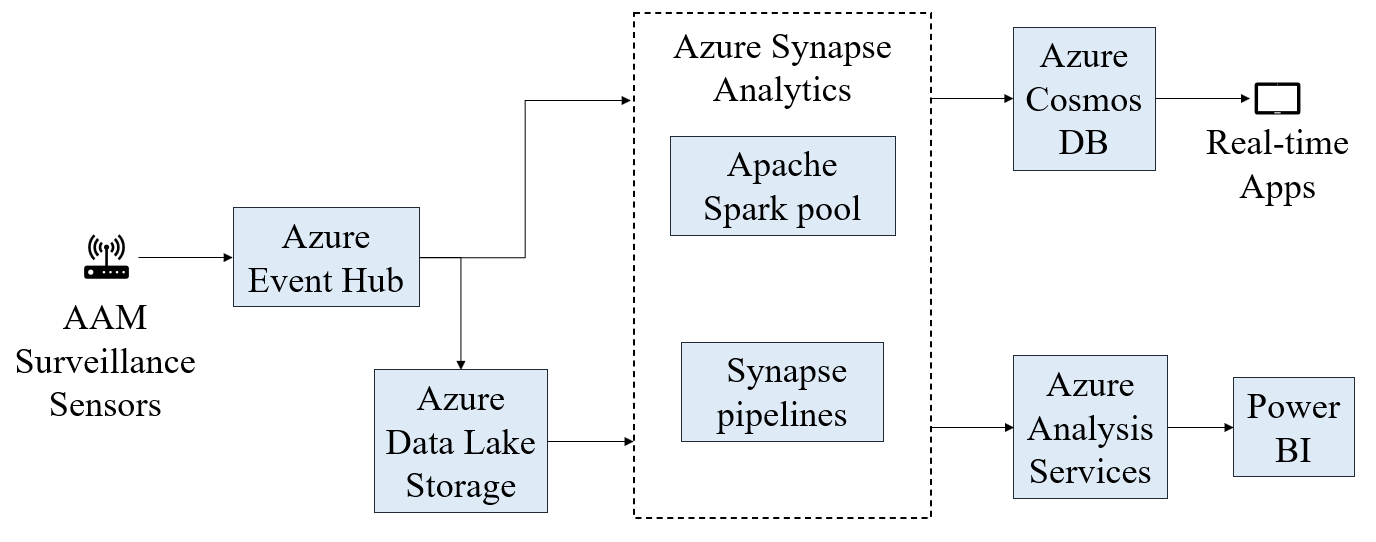}
\end{adjustwidth}
    \caption{A flowchart showing the connections of the cloud components of LASIC.}
    \label{fig:cloud}
\end{figure}

A cloud computing architecture capable of real time analytics on big data would need to be created to enable the data flow through LASIC. The cloud computing architecture would consist of six components: 1) Azure Event Hub, 2) Azure Synapse Analytics, 3) Azure Data Lake Storage, 4) Azure Cosmos Database (DB), 5) Azure Analysis Services, and 6) Power BI \cite{azure}. The cloud components of LASIC are depicted in Figure \ref{fig:cloud}. The Azure Event Hub is a big data streaming platform and event ingestion service, where millions of data units can be received and processed in a single second \cite{eventhub}. It can be used to easily ingest live streaming data from the AAM surveillance sensors. Then, Azure Synapse Analytic can be used to transform and store data that has been provided to the Azure Event Hub, respectively. The Azure Synapse Analytics is an analytics service that combines data integration, enterprise data warehousing, and big data analytics \cite{synapse}. For large-scale access and movement of surveillance data, Azure Synapse Analytics would require the use of Apache Spark pool and Synapse pipelines. These components can be used for data cleaning, transforming, and analyzing; and can enable the use of Python, Scala, or .NET, and scalable ML techniques to derive deeper insights from LASIC data. Azure Data Lake Storage allows massively scalable and secure data lake functionality built on Azure Blob Storage \cite{storage}, which is needed to store the LASIC data. To provide access to the intended LASIC data to subscribers through real-time apps, data would need to be transferred from Apache Spark pools to Azure Cosmos DB \cite{cosmos}. Analytics dashboards and embedded reports on LASIC data can be created using Azure Analysis Services and Power BI for use by the LASIC operator and subscribers \cite{analysis, power}. 

% The Azure Blob Storage helps to create data lakes for analytics needs, and provides data storage \cite{blob}.

\subsection{Cost-Benefit Analysis of Low Altitude Surveillance Information Clearinghouse}\label{CBA}

To assess and justify the worth of investing in an AAM surveillance network and LASIC, a cost-benefit analysis is needed, as mentioned in Section \ref{intro}. Hence, we conduct the cost-benefit analysis which involves identifying and estimating the potential costs and benefits associated with this infrastructure. The findings of the first two parts are used to perform the analysis, as shown in Figure \ref{fig:flow}. The cost and benefit factors identified and considered to be significant in this analysis are 1) surveillance sensor cost (cost factor 1), 2) cloud computing cost (cost factor 2) and 3) revenue generated from LASIC subscription (benefit factor), which are discussed more in this section. Using the estimates of the cost and benefit factors, NPV is calculated, which is necessary to measure of the future return on investment expected from an investment in a project in terms of today's dollars. The NPV metric takes into account the time value of money and future cash flows, which is further discussed in Section \ref{NPV}. As there is uncertainty associated with some of the key parameters involved in the analysis --- namely, subscription fees, number of initial subscribers, terrain-based sensor detection probabilities and minimum required detection probability --- a sensitivity analysis is performed to evaluate the effect of these parameters on the NPV generated.

 % These inputs along with AAM traffic projections, cloud computing policies and the output from SAND model were used to estimate the cost and benefit factors.
% and LASIC functional and performance requirements
% incorporated in the computation of the cost and benefit factors, as discussed below.

% Next, we estimated the cloud computing cost using the AAM traffic projections data and Microsoft Azure cloud computing pricing policies. 

% From SAND model, we determined where and how many sensors of each sensor type would be required in SMCO to build AAM surveillance network. 

% Based on these features and functionalities, the cost and benefit factors were identified and evaluated and then  

\subsubsection{Cost Factor 1: Surveillance Sensor Cost} 

The first cost factor is the surveillance sensor cost, which is total cost incurred to build the surveillance network in a given area. This cost factor depends on the number of sensors needed to get the intended coverage over a given region and the price of sensors. The required number of sensors is evaluated by implementing SAND model described in Section \ref{model} for six sensor types mentioned in Section \ref{sensortype}. In this analysis, we consider the capital for purchasing and installing sensors to build the surveillance network will be invested once in 2023.

\subsubsection{Cost Factor 2: Cloud Computing Cost} \label{cf2}

As discussed in Section \ref{ccc}, opting to host LASIC on a cloud server incurs an annual cloud computing cost that constitutes the second cost factor throughout the analysis period. To estimate this cost, we refer to Microsoft Azure's pricing policies, as outlined in \cite{azure_price}, which take into account the number of surveillance data published and received by LASIC. The amount of surveillance data generated in a particular area can vary depending on factors such as AAM traffic projection, surveillance data type, interface standard, data size and ping rate.

\subsubsection{Benefit Factor} \label{benefit_fact}

The survey responses reveal the willingness to pay of potential subscribers of LASIC for the services offered by it. The range of subscription fees that potential subscribers are willing to pay for the services expected to be offered by LASIC is found to be \$100-\$400. This inform the computation of the benefit factor considered in this study, the revenue generated from LASIC. The potential subscribers of LASIC include parcel and cargo delivery operators, medical item delivery companies, air taxi operators, infrastructure inspection companies, airspace service providers, state penitentiaries, law enforcement agencies, correctional facilities, and municipalities. The number of potential subscribers in the various years of the analysis period are estimated based on global and US AAM market growth rates reported in AAM market studies such as \cite{cagr}, \cite{cagr1}, and \cite{cagr2}. 

% The benefit factor considered in the cost-benefit analysis of AAM surveillance network and LASIC is the revenue generated from subscription fees charged to the potential subscribers of LASIC. According to the survey responses, 

\subsubsection{Net Present Value} \label{NPV}

% To determine if a project will generate enough revenue to cover its expenses, NPV is a highly relevant metric. NPV expresses the total amount of money that an entity can expect to earn from an investment over its lifetime, including both future positive and negative cash flows \cite{dobrowolski2022does}. If NPV is positive, it indicates that the projected revenues from the project or investment will be greater than the estimated costs in current dollars. An investment with a positive NPV is considered profitable, while one with a negative value is considered to result in a net loss. The projected yearly NPV depends on yearly cash inflow and cash outflow, where cash inflow refers to revenue generated from subscription fees and cash outflow refers to the costs of surveillance sensors and cloud computing \cite{altonji2023effects}. The NPV calculation uses a discount factor to account for the time value of money, which reflects the idea that a dollar received in the future is worth less than a dollar received today. It means that the value of money decreases over time due to the effects of inflation, opportunity cost, and other factors. 

To evaluate the financial viability of LASIC, its NPV over the analysis period needs to be computed and analyzed. 

\begin{equation} \label{11}
NPV_t =  \frac{C^P_t-C^N_t}{{(1+\chi)^{t}}},  \indent \forall t \in \{0, 1, ..., 10\} 
\end{equation}

\noindent represents the estimated total value of all future cash flows generated by an investment over the lifetime of the project or a given analysis period, where $C^P_t$ and $C^N_t$ represent positive cash flow and negative cash flow in year $t$, respectively \cite{gallo2014refresher, dobrowolski2022does, altonji2023effects}. In this analysis, the yearly NPV calculation of AAM surveillance network and LASIC is carried out based on the difference between the revenue generated (positive cash flow) and the costs associated with AAM surveillance network and LASIC (negative cash flow). A discount factor is considered to account for the time value of money, which reflects the idea that a dollar received in the future is worth less than a dollar received today. We consider cash flow over a 10 year horizon, discounting at 10\% (the discount rate $\chi$) for this infrastructure project \cite{terrill2018unfreezing, jawad2006discount, nokkala2004role}. A positive NPV at the end of the analysis period implies that the expected revenue from the investment exceeds the projected costs, and thus, the investment is considered profitable. Conversely, a negative NPV suggests that the investment would result in a net loss.

\section{Results}\label{sec4}

% The SAND model was implemented in Python using the Gurobi-Python API and the model was solved using the Gurobi Optimizer. Based on the solutions of the SAND model, this section was designed into four sub-sections to offer a comprehensive analysis of the framework of AAM surveillance sensor network and LASIC. The first sub-section is the revenue analysis, which focuses on evaluating the financial performance of the sensor network with regards to its revenue. The second sub-section, homogeneous sensor analysis, presents the examination of each individual sensor type used in the study. Surveillance sensor cost and cloud computing cost were calculated to consider the cost factors. To analyze the insights of each sensor type, the NPV plots were compared with each other in terms of two criteria: 1) the number of years required to reach the break-even point (BEP), where NPV is equal to zero, and 2) the estimated NPV in the final year of the analysis period. The third sub-section, heterogeneous sensor analysis, analyzes the behavior of the surveillance sensor network when different types of sensors are combined. Lastly, the sensitivity analysis examines how the framework responds to changes in key parameters and assesses the sensitivity of the framework to these changes.

The applicability of the proposed model and solution approaches are demonstrated through their application in several numerical experiments in this section. We considered two distinct types of surveillance sensor networks for the numerical experiments, namely homogeneous and heterogeneous sensor networks. Homogeneous sensor network comprises only one type of sensor whereas heterogeneous sensor network is composed of different types of sensors. To undertake the experiments, we consider the six major cities of Ohio (SMCO): Columbus, Cleveland, Cincinnati, Akron, Toledo, and Dayton. The formation of SMCO is predicated on the finding of significant demand potential for AAM use cases in those cities considering socioeconomic factors, such as population, population density, gross domestic product, median per capita income, cost of living, total area, cities in motion index, human capital, etc \cite{del2021infrastructure}.

% \subsection{Numerical Experiments}
\subsection{Experimental Setup}

In this section, we provide details of the values related to sensors and surveillance area that we consider in our experimental setup for running the experiments.

\subsubsection{Sensors}

In this study, one real-world sensor model is considered for each of the sensor types discussed in Section \ref{sensortype}. The Echo Guard radar is considered for the radar sensor type. It is a top-tier 4D radar with an easy user interface that is easily adaptable to site and mission requirements for high performance ground-based detect and avoid \cite{echo}. We consider CamelCase pingStation3 from AVIONIX Software S.L. as an ADS-B frequency ground receiver for our analysis. It is a networkable weatherproof 978/1090 MHz ADS-B receiver including GPS and antenna, with power and data provided by a single power-over-Ethernet network cable connected right to LAN \cite{ads-b}. We consider DroneScout as the (direct/broadcast) remote ID receiver, which can receive remote ID signals sent from aircraft \cite{remote}. The Dedrone RF-360 is considered in this study for the RF sensor type. It is a passive, network-attached radio sensor for the detection, classification, and localization (geolocation) of aircraft and their remote controls \cite{RF360}. The OptiNav Drone Hound system is considered as the acoustic sensor type in our study. It is an acoustic sensor that has been designed to detect, identify, and track sUAS. Unlike other sensors, it does not rely on electromagnetic emissions from the sUAS \cite{aco}. We consider the Q6225-LE PTZ Network Camera from Axis Communications in our analysis for the optical camera sensor type \cite{camera}. 

% In this study, one real-world sensor model is considered for each of the sensor types discussed in Section \ref{sensortype}. The Echo Guard radar is considered for the radar sensor type. It is a top-tier 4D radar with an easy user interface that is easily adaptable to site and mission requirements for high performance ground-based detect and avoid (see Figure \ref{fig:comb}a). We consider Ping Station 3 from AVIONIX Software S.L. as an ADS-B frequency ground receiver for our analysis (see Figure \ref{fig:comb}b). It is a networkable weatherproof 978/1090 MHz ADS-B receiver including GPS and antenna, with power and data provided by a single power-over-Ethernet network cable connected right to LAN \cite{ads-b}. We consider DroneScout as the (direct/broadcast) remote ID receiver, which can receive remote ID signals sent from aircraft (see Figure \ref{fig:comb}c) \cite{remote}. The Dedrone RF-360 is considered in this study for the RF sensor type. It is a passive, network-attached radio sensor for the detection, classification, and localization (geolocation) of aircraft and their remote controls (see Figure \ref{fig:comb}d) \cite{RF360}. The OptiNav Drone Hound system is an acoustic sensor that has been designed to detect, identify, and track sUAS. Unlike other sensors, it does not rely on electromagnetic emissions from the sUAS (see Figure \ref{fig:comb}e). We consider the Q6225-LE PTZ Network Camera from Axis Communications in our analysis for the optical camera sensor type (see Figure \ref{fig:comb}f) \cite{camera}. 

The $R$, $\psi$, and $\delta$ of the six selected sensors from their corresponding sensor types in the input set $S$ of the SAND model are listed in Table \ref{tab:sensor}. The sensor types listed vary in terms of their range and cost. Radar has the highest range of 321.87 km, provided by the ADS-B sensor. Remote ID and RF sensors have ranges of 5.02 km and 4.99 km respectively, while acoustic and optical camera sensors have much shorter ranges of 0.5 km and 0.4 km, respectively. Radar and RF sensors are generally the most expensive, with one of the sensors costing around \$35,000. The other sensors in the table range in price from \$1,100 to \$9,000.

\begin{table}[H]
\caption{\label{tab:sensor} Selected sensors from each sensor type.}
\begin{adjustwidth}{-\extralength}{0cm}
		\newcolumntype{C}{>{\centering\arraybackslash}X}
		\setlength{\tabcolsep}{4.8mm}
		\begin{tabularx}{\fulllength}{lcccccc}
\toprule
% & Transition& & \multicolumn{2}{c}{}\\\cmidrule{2-2}
\textbf{Sensor Types} & \textbf{Vendor} & \textbf{System} & \textbf{Range} &  $\boldsymbol{\psi}$ & $\boldsymbol{\delta}$ \\
  & & & \textbf{(km)} & \textbf{($\approx$ \$)}&  \\
\midrule
Radar & Echodyne & Echo Guard & 2.41& 35000 & 3\\
ADS-B & AVIONIX Software S.L.  & CamelCase pingStation3 & 321.87& 2250 & 1\\
Remote ID & BlueMark Innovations BV & Drone Scout & 5.02& 1100 & 1\\
RF & Dedrone & RF - 360 & 4.99& 35000 & 1\\
Acoustic & OptiNav & Drone Hound & 0.5& 9000 & 1\\
Optical Camera &  Axis Communications &  Q6225-LE PTZ Network Camera & 0.4& 3500 & 6\\
\bottomrule
\end{tabularx}
\end{adjustwidth}
\end{table}

\subsubsection{Surveillance Area}

In the experimental setup for modeling a surveillance area, we consider the following factors to execute the SAND model. As mentioned in Section \ref{model}, $L/\sqrt{2}$ should be greater than $\rho$. Therefore, we consider $L$ as 0.3 km, which is less than 0.4$\sqrt{2}$ km, where 0.4 km is the lowest value of range among all sensor types $S$ (refer to Table \ref{tab:sensor}). This ensures that all blocks are covered by the sensors and no surveillance area is left uncovered. Hence, to run SAND model for Columbus, Cleveland, Cincinnati, Toledo, Akron, and Dayton, we generate 130 × 126,  77 × 96, 58 × 77, 57 × 72, 65 × 62, and 54 × 60 blocks, respectively.

% We considered a value of 0.3 km for $L$ when executing the SAND model in our study. This value is less than 0.4$\sqrt{2}$ km, where 0.4 km represents the lowest value of $R$ among all sensor types $S$ .

A higher minimum required detection probability value would need to be considered for security-sensitive areas to ensure a higher level of detection and identification accuracy. We consider the value of $r$ as 0.98 in our study, which means that the system aims to detect and identify targets with a probability of at least 0.98. For considering the effect of terrains on the probability of detection of a sensor, the terrains of the area are divided into five major types of $T$: 1) open, 2) water, 3) neighborhood or residential area, 4) hill, and 5) busy commercial area or downtown. The Google Maps platform is utilized to observe and determine the terrain type of each block $z \in Z$. $T_z$ is obtained by reference to the terrain classification of the $z$-th block, as recorded in the list $T$. The probability values in the $\omega^s_T$ matrix, as given in Table \ref{tab2}, are approximated based on the approach reported in \cite{lamm2001develop, seo2008optimal, seo2016efficient}. By analyzing the values presented in the table, it is evident that the detection probabilities of the sensors tend to decrease as the terrain type changes, which aligns with the discussion in Section \ref{detect_prob}. The probability of detection is the highest for open terrains, followed by water and neighborhood, whereas it is comparatively lower for hill and commercial areas. Additionally, the values indicate that the detection probabilities for radar and ADS-B sensor are relatively higher than other sensor types, regardless of the terrain type. Conversely, acoustic sensor has the lowest detection probability among all sensor types and terrain types. These observations emphasize the significance of selecting an appropriate sensor type and its detection probability when designing a sensor network for a given area with varying terrain types.

\begin{table}[H]
\caption{\label{tab2} Detection probability of sensors based on different terrain types: $\omega^s_T$ matrix.}
\newcolumntype{C}{>{\centering\arraybackslash}X}
\setlength{\tabcolsep}{3.3mm}
\begin{tabularx}{\textwidth}{lcccccc}
\toprule
&&&&Terrain Type ($T$)& \multicolumn{2}{c}{}\\\cmidrule{2-7}
Sensor Type ($S$) && Open & Water& Neighborhood & Hill & Commercial Area\\\midrule
Radar && 0.95 & 0.90 & 0.85 & 0.75 & 0.75\\
ADS-B && 0.99 & 0.99 & 0.90 & 0.85 & 0.80\\
Remote ID && 0.95 & 0.95 & 0.85 & 0.80 & 0.75\\
Radio Frequency && 0.95 & 0.95 & 0.85 & 0.80 & 0.75\\
Acoustic && 0.75 & 0.65 & 0.40 & 0.25 & 0.20\\
Optical Camera && 0.90 & 0.90 & 0.80 & 0.75 & 0.70\\
\bottomrule
\end{tabularx}
\end{table}

% A heatmap showing the probability of detection of a radar across the different terrain types in Dayton is presented in Figure \ref{fig4}. 

\begin{figure}[H]
    \includegraphics[width=10.5cm,height=8cm]{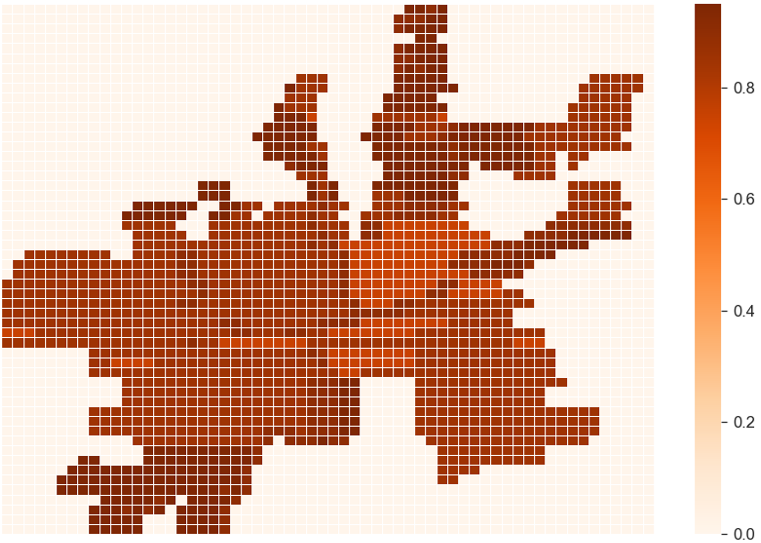}
    \caption{A heatmap of probability of detection of a radar based on terrain types of 3240 blocks in Dayton.}
    \label{fig4}
\end{figure}

To demonstrate the dependence of a sensor's detection probability on the terrain type, Figure \ref{fig4} presents a heatmap showing the probability of detecting radar across the various terrain types in Dayton. The colored bar on the right side of the figure shows the scale of detection's probability, where the off-white color refers to a probability of zero and the darkest orange color refers to a probability of one. Note that the detection probability of sensors in blocks outside the area, where $I(z) = 0$, is set to zero. Hence, the off-white color refers blocks that do not belong to Dayton, and the orange colors, on the other hand, represent blocks within the area. By setting the detection probability of sensors in blocks outside the area to zero, the analysis is focused on the sensors within the area of interest. This allows for a more precise evaluation of the sensor network's effectiveness in the designated area.

\begin{figure}[H]
\includegraphics[width=10cm,height=7cm]{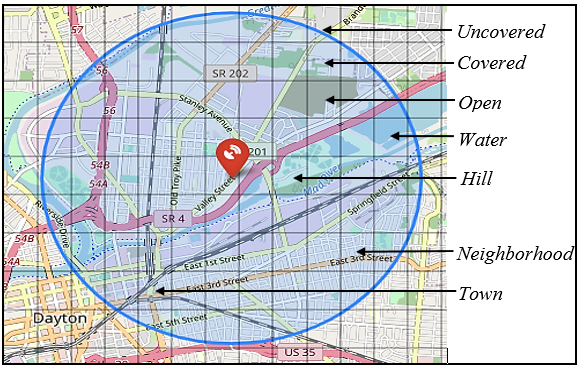}
\caption{Block selection overview - classification of blocks by sensor coverage and terrain type.}
\label{fig5}
\end{figure}

%updated version
Considering the sensor detection probabilities in different terrain types, Figure \ref{fig5} presents an overview of the selection process that the SAND model uses to select blocks within a sensor range. The red marker refers to the location where a sensor is placed, and the blue circle shows the area within its range. Blocks which do not have all four corner points inside the blue circle are classified as the "Uncovered" blocks (i.e., uncovered by the sensor range), whereas the blocks with all four corner points inside the blue circle are classified as the "Covered" blocks. Each of the covered blocks can belong to any of the five terrain types in $T$, as mentioned in Figure \ref{fig5} in italic font.

%old version (please see the updated version above)
% An overview is given in Figure \ref{fig5} to depict how SAND model selects the blocks within a sensor range. The red marker refers to the location where a sensor is placed, and the blue circle shows the area within its range. Blocks which do not have all four corner points inside the blue circle are considered to be uncovered blocks (i.e., uncovered by the sensor range), whereas the blocks with all four corner points inside the blue circle are considered as covered blocks. There can be five types of covered blocks according to the terrain types which were mentioned in Figure \ref{fig5} in italic font. 

Following the methodology described in Section \ref{CBA}, determining the number of sensors needed and their locations using the SAND model, the revenue, costs and NPV of AAM surveillance network and LASIC are computed for the numerical experiments. Before discussing the homogeneous and heterogeneous sensor placement analysis, we address the analysis of revenue and cloud computing cost since they are the same for both experiments. NPV of the different sensor types are compared with each other in terms of two criteria: 1) the number of years required to reach the break-even point, and 2) the estimated NPV in the final year of the analysis period. Given the uncertainty of AAM market, a sensitivity analysis is conducted to examines how the NPV responds to changes in key market parameters --- namely, the yearly number of subscribers and subscription fees. The results are presented in the following subsections.

\subsection{Revenue and Cloud Computing Cost Analysis} \label{revenue_cloud}

The total yearly revenues generated by AAM surveillance network and LASIC is determined by two factors: the yearly number of subscribers and the monthly subscription fee of LASIC. It does not depend on the sensor type used in the surveillance network provided that complete coverage is present across SMCO. Hence, for all sensor types, the revenue generated is the same. We estimate the number of subscribers for LASIC in 2024 to be 100 based on the existing number of sUAS operators, air taxi operators, infrastructure inspection companies, medical item delivery companies, state penitentiaries, law enforcement agencies, correctional facilities, and municipalities in SMCO. Then, we study the compound annual growth rate (CAGR) of the market size of AAM to estimate the number of yearly subscribers of LASIC over the analysis period based on the reports available on global and US AAM market growth \cite{cagr, cagr1, cagr2, cagr3}. Because of the evolving nature of the AAM market and its inherent uncertainty, we incorporate a CAGR range of 10\% to 20\% in our study based on the values reported in the AAM market growth studies instead of relying on a fixed value of CAGR. Here, the lower limit of 10\% signifies the conservative estimate, while the upper limit of 20\% represents the optimistic estimate.

% Then, we follow a 15\% compound annual growth rate (CAGR) of the market size of AAM to estimate the number of yearly subscribers of LASIC over the analysis period based on the reports available on global and US AAM market growth \cite{cagr, cagr1, cagr2}. To consider the uncertainty in the AAM market, we incorporate a CAGR range of 10\% to 20\% in our study. 

% Then, we followed a 15\% compound annual growth rate (CAGR) of the market of eVTOL for passenger and heavy cargo transportation, and a 45\% CAGR of the market of sUAS for other use cases to estimate the number of yearly subscribers of LASIC over the analysis period based on the reports available on global and US AAM market growth  \cite{cagr, cagr1, cagr2}.

\begin{figure}[H]
{\captionsetup{position=bottom,justification=centering}
\begin{subfigure}{.48\textwidth}
  \centering
  % include first image
  \includegraphics[width=1\linewidth]{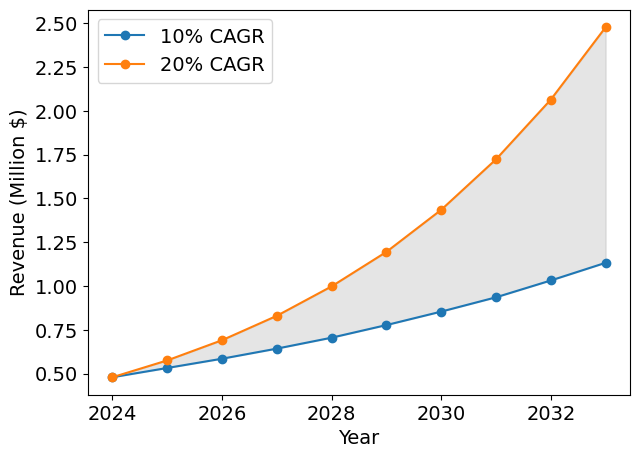}
  \caption{Yearly revenues generated by AAM surveillance network and LASIC.}
  \label{fig:revenue}
\end{subfigure}
\begin{subfigure}{.49\textwidth}
  \centering
  % include second image
  \includegraphics[width=1\linewidth]{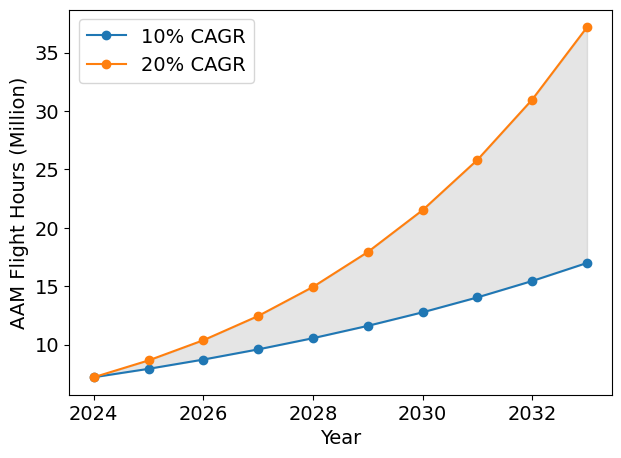}  
  \caption{Yearly projected AAM traffic in Ohio.}
  \label{fig:AAM_flight}
\end{subfigure}}
\caption{Yearly revenues and yearly projected AAM traffic.}
\label{fig:rev_fli}
\end{figure}

Another factor that affects the revenue is the monthly subscription fee, which ranges between \$100 and \$400, as mentioned in Section \ref{benefit_fact}. For the revenue analysis, we assume the fee to be \$400 and later vary it during sensitivity analysis to analyze cases where the fee is less than \$400. The yearly revenues generated by AAM surveillance network and LASIC in SMCO, with a initial number of subscribers of 100 and fixed subscription fee of \$400, are depicted in Figure \ref{fig:revenue}. As the number of subscribers increases over the years, the revenue grows proportionally. The shaded region signifies the potential revenue outcome that falls between the projected revenues at a 10\% CAGR and those at a 20\% CAGR. It also highlights how the growth rate significantly impacts the range of revenue projections over time. In the case of a 10\% CAGR, the revenue is expected to begin at \$0.48 million and increase gradually over the years, reaching approximately \$1.032 million in 2033. On the other hand, with a more optimistic 20\% CAGR, the revenue starts at the same initial value of \$0.48 million but experiences a more rapid growth, reaching a potential high of \$2.064 million by 2033.

% The revenue starts at \$0.6 million in the first year and is projected to reach approximately \$15 million by 2033.

% revenue plot
% \begin{figure}[H]
%     \centering
%     \includegraphics[width=9cm,height=6.5cm]{revenue.PNG}
%     \caption{Yearly revenues generated by AAM surveillance network and LASIC.}
%     \label{fig:revenue}
% \end{figure}

\begin{figure}[H]
{\captionsetup{position=bottom,justification=centering}
\begin{subfigure}{.48\textwidth}
  \centering
  % include first image
  \includegraphics[width=1\linewidth]{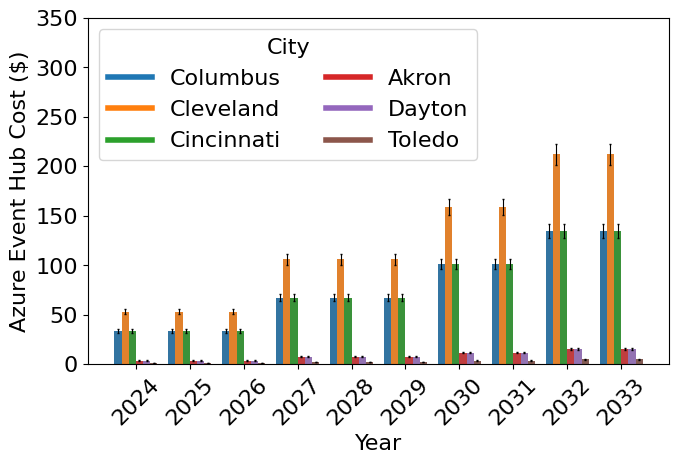}
  \caption{Azure Event Hub cost}
  \label{fig8a}
\end{subfigure}\vspace{8pt}
\begin{subfigure}{.48\textwidth}
  \centering
  % include second image
  \includegraphics[width=1\linewidth]{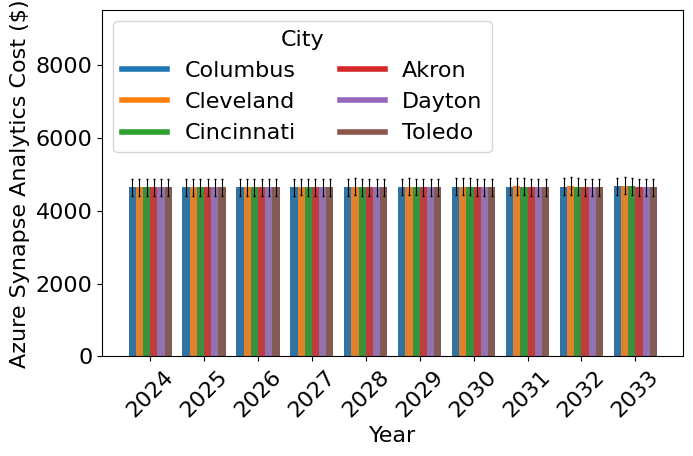}  
  \caption{Azure Synapse Analytics cost}
  \label{fig8b}
\end{subfigure}
\vspace{8pt}
\begin{subfigure}{.48\textwidth}
  \centering
  % include third image
  \includegraphics[width=1\linewidth]{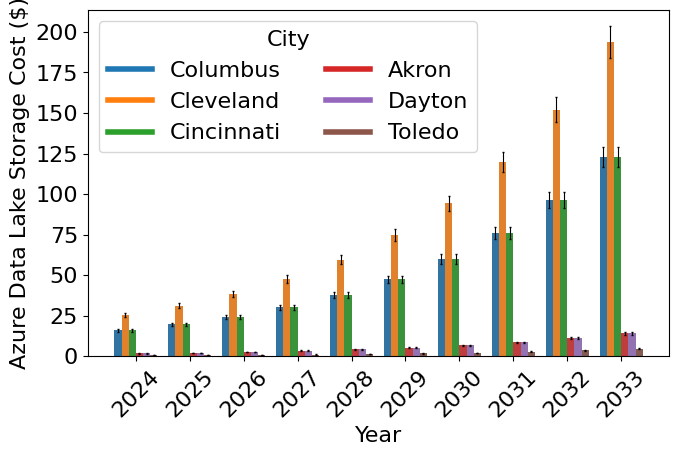}  
  \caption{Azure Data Lake Storage cost}
  \label{fig8c}
\end{subfigure}
\begin{subfigure}{.48\textwidth}
  \centering
  % include fourth image
  \includegraphics[width=1\linewidth]{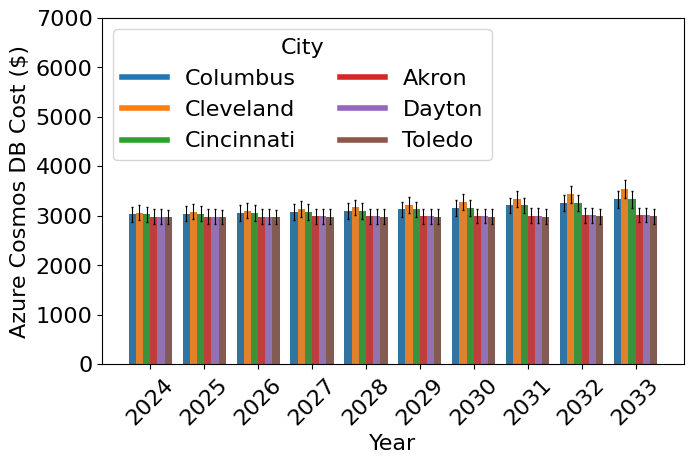}  
  \caption{Azure Cosmos DB cost}
  \label{fig8d}
\end{subfigure}}

\centering
\begin{subfigure}{.48\textwidth}
  \centering
  % include third image
  \includegraphics[width=1\linewidth]{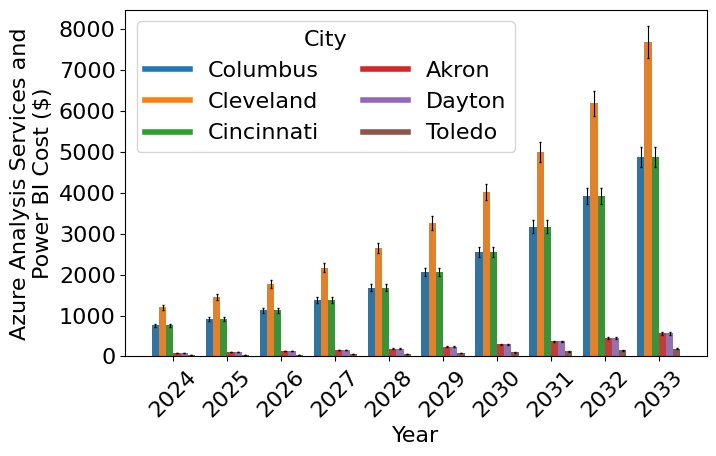}  
{\captionsetup{position=bottom,justification=centering}\caption{Azure Analysis Services and Azure Power BI cost}
  \label{fig8e}}
\end{subfigure}
\caption{Cost of different cloud components in cloud computing.}
\label{fig8}
\end{figure}

% We consider a range of 10\% to 20\% in the CAGR of AAM market size to account for uncertainty in AAM flight hours
% should fall
The cost of cloud computing, according to the Microsoft Azure pricing policies, depends on the projected amount of surveillance data generated in each city, which in turn depends on the projected AAM traffic in each city as presented in Figure \ref{fig:AAM_flight}. The range of  CAGR values of AAM market size considered accounts for the uncertainty in AAM flight hours. This uncertainty range subsequently influences the cloud computing cost, as illustrated in the figures through the inclusion of error bars. These error bars serve to indicate that the associated cost is estimated to fall within the limits defined by the bar. The ten-year cloud computing cost breakdown for each cloud component and each city within SMCO is illustrated in Figure \ref{fig8}. The Azure Event Hub and Azure Data Lake Storage have the two lowest costs among all the components. The Azure Event Hub operates on a tiered pricing model, where the cost of the service varies based on the level of usage of surveillance data by a subscriber. The cost is relatively low in the first few years as the surveillance data and number of subscribers is initially low, and the cost increases in steps as as the surveillance data and number of subscribers increases. When the Azure Event Hub usage reaches a defined threshold, the cost climbs to a higher level, and this pattern repeats for each subsequent tier, creating a step function of the cost with respect to usage. The costs for Azure Event Hub and Azure Data Lake Storage increases with the amount of incoming data ingested into the hub and stored in the Data Lake, respectively. Additionally, the frequency of data access also influences the rise in cost, with higher amounts of access due to increasing number of subscribers leading to an increase in cost in successive years. The Azure Analysis Services and Azure Power BI costs increase with time commensurate with the projected increase in the amount of data stored, number of queries run, and number of users accessing the services. Lastly, the pricing of Azure Synapse Analytics and Azure Cosmos DB includes both a yearly fixed cost and a yearly variable cost. The yearly fixed cost is associated with the provisioning of virtual machines, storage, and other necessary resources to operate the services. The yearly variable cost depends on the amount of data processed in LASIC. As the yearly fixed cost is much higher than the yearly variable cost, Azure Synapse Analytics and Azure Cosmos DB costs are nearly constant, increasing slightly over the years. Across all cities, the cloud computing cost associated with ingesting, storing, and analyzing the surveillance data generated in Cleveland is the highest as it has the highest air traffic demand forecast across SMCO, and hence produces the the largest amount of surveillance data; whereas for Toledo, the cost is the lowest as it generates the lowest air traffic demand forecast, and hence the lowest amount of data. 

% cloud computing cost plots

\subsection{Homogeneous Sensor Placement Analysis} \label{homo}

In the homogeneous sensor placement analysis, the surveillance network across SMCO is considered to be built using one sensor type instead of a mix of sensor types. This allows for a more in-depth analysis of each individual sensor type's suitability for AAM surveillance and capability to produce NPV over the analysis period. 

\begin{figure}[H]
\begin{adjustwidth}{-\extralength}{0cm}
%\centering %% If there is a figure in wide page, please release command \centering
    \centering
   \includegraphics[width=15cm,height=9cm]{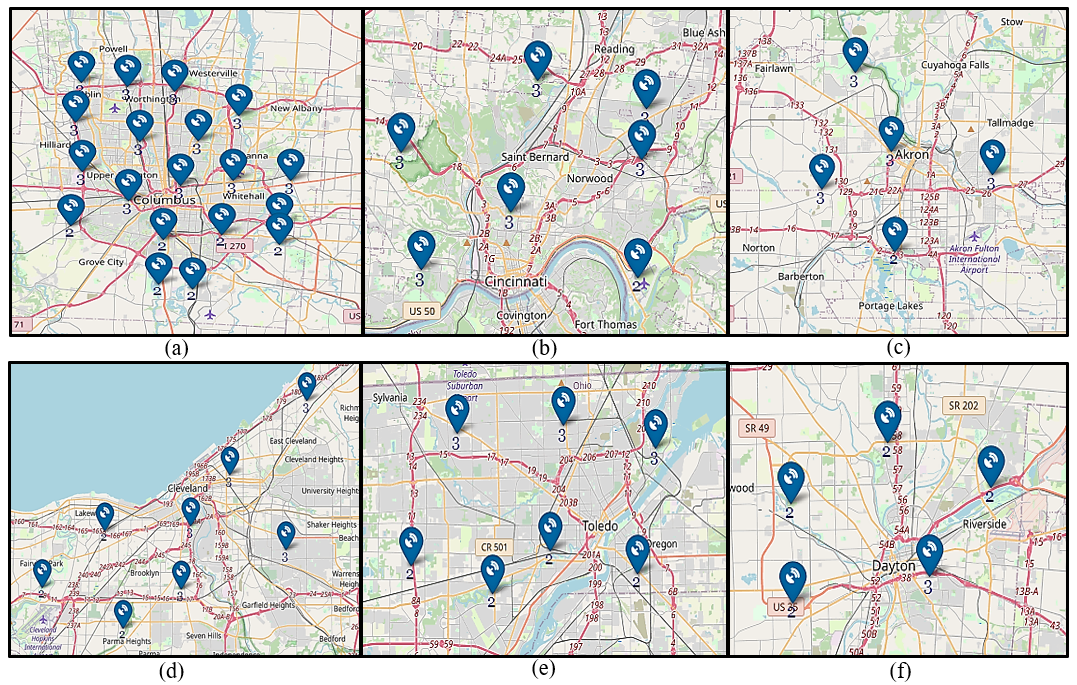} 
\end{adjustwidth}
    
    \caption{Optimal locations of RF sensors (the blue markers) in six cities (a) Columbus, (b) Cincinnati, (c) Akron, (d) Cleveland, (e) Toledo, and (f) Dayton.}
    \label{fig6}
\end{figure}

For each sensor type, the optimal location and number of sensors required to build the homogeneous surveillance network at minimum sensor cost in SMCO are determined using the SAND model. The optimal locations for RF sensors in the surveillance network across SMCO are shown in Figure \ref{fig6}, where the blue markers represent the locations of the RF sensors.
% , and the text below the markers indicates the number of RF sensors needed at each location. 
The sensors are strategically placed to ensure that all blocks within the sensor range are covered while minimizing the overlapping region to reduce the cost of the sensor network. The distribution of sensor locations varies for each city based on factors such as city shape and area, terrain type, and probability of detection of each sensor type based on terrain types. These factors also affect the optimal location and number of sensors required to provide adequate coverage of the city and determine the total sensor cost. For example, Columbus can be approximated as having a circular shape, while Cleveland is wider than it is long, and Columbus has a larger area compared to cities like Akron and Dayton. Another instance is Toledo, which has a greater number of water bodies compared to Columbus, while Cincinnati has more hilly terrain. These variations in shape and terrain result in different sensor distribution in the network.

 % (the text below blue marker shows the number of sensors needed at each location)

% Based on the number of sensors required for each city as determined by the SAND model, we calculated the city-wise sensor cost of six sensor types as shown in Figure \ref{fig7}. Though the price of one acoustic sensor is cheap, it requires a large number of sensors to cover SMCO because of its small range. Hence, it becomes the most expensive sensor type in terms of total sensor cost, while ADS-B is the least expensive among the six types. The sensor types are listed in ascending order of cost as follows: ADS-B, Remote ID, RF, Radar, Optical Camera, and Acoustic. Among the six cities, Columbus requires the highest sensor cost and Dayton the lowest.

% sensor locations plot

The number of sensors required varies significantly depending on the city and the type of sensor used, as shown in Table \ref{city_sensors}. For any given sensor type, the number of sensors required to cover a given area increases with the area of the city. Among SMCO, Columbus requires the largest number of sensors as it has the largest area and Dayton the smallest as it has the smallest area. In addition, the range of the sensors used also affects the number of sensors required. Radar sensors, for example, typically have a longer range than optical cameras (refer to Table \ref{tab:sensor}), which means that fewer radar sensors are needed to cover the same city compared to using more optical cameras. Different cities have varying proportions of terrain types (refer to Figure \ref{fig6}). The detection probabilities of sensors on different terrain types are affected by the sensor type used (refer to Table \ref{tab2}), which in turn affects the number of sensors required. For instance, acoustic sensors have a lower detection probability than other sensor types when placed on a block with hilly terrain. As a result, cities with hilly terrain, such as Cincinnati, require more sensors to cover the terrain than cities like Cleveland with relatively less hilly terrain. Therefore, the number of sensors needed in Cincinnati is higher than the number of sensors needed in Cleveland for acoustic sensors, even though Cleveland is larger in size compared to Cincinnati. Moreover, the field of view varies with the sensor type, and for a limited field of view, more sensors are required to ensure a complete \SI{360}{\degree} view. For instance, although the ranges of acoustic sensors and optical cameras are similar, acoustic sensors have a lower detection probability compared to optical cameras. Thus, the number of sensors required for acoustic sensors to cover an area should be higher than that required for optical cameras. However, due to the higher value of $\delta$ for optical cameras compared to acoustic sensors (refer to Table \ref{tab:sensor}), the number of sensors required becomes higher for optical cameras than for acoustic sensors.

\begin{figure}[H]
{\captionsetup{position=bottom,justification=centering}
\begin{subfigure}{.48\textwidth}
  \centering
  % include first image
  \includegraphics[width=1\linewidth]{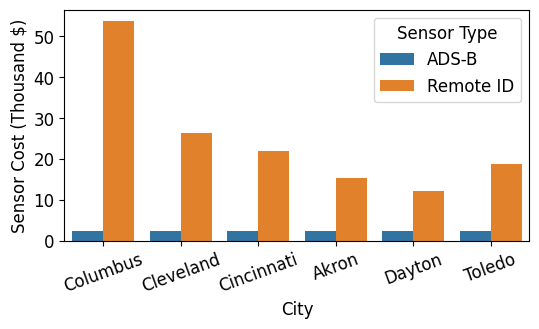}
  \caption{Sensor cost of ADS-B and remote ID.}
  \label{fig7a}
\end{subfigure}\vspace{8pt}
\begin{subfigure}{.48\textwidth}
  \centering
  % include second image
  \includegraphics[width=1\linewidth]{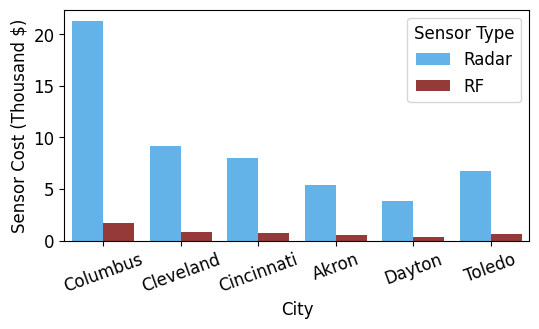}  
  \caption{Sensor cost of radar and RF.}
  \label{fig7b}
\end{subfigure}}

\centering
\begin{subfigure}{.48\textwidth}
  \centering
  % include third image
  \includegraphics[width=1\linewidth]{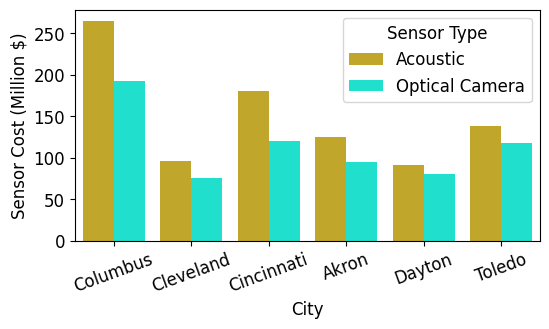}  
  \caption{Sensor cost of acoustic and optical camera.}
  \label{fig7c}
\end{subfigure}
\caption{City-wise sensor cost for different sensor types.}
\label{fig7}
\end{figure}

Based on the unit price of each sensor type and number of sensors required for a sensor type for each city, the SAND model determines the city-wise sensor cost of all sensor types. The sensor cost of all cities are presented in Figure \ref{fig7}. The different sensor types listed in ascending order of sensor cost are: ADS-B, remote ID, RF, radar, optical camera, and acoustic. Among SMCO, Columbus requires the highest sensor cost for all sensor types, while Dayton has the lowest. The sensor costs for ADS-B and remote ID sensor types are much lower compared to the other sensor types, as they require fewer sensors and have lower unit prices. At the other end of the sensor cost spectrum are optical cameras and acoustic sensor types. Though the unit prices of acoustic sensors and optical cameras are cheap, a large number of sensors are required for each to cover SMCO, as mentioned above, resulting in a very high total sensor cost.  While the number of sensors required for the acoustic sensor type is less than for optical cameras, the total sensor cost for acoustic sensors is higher than for optical cameras due to the lower unit price of optical cameras.

\begin{table}[H]
\caption{Number of sensors required in SMCO.}
\label{city_sensors}
\newcolumntype{C}{>{\centering\arraybackslash}X}
\begin{tabularx}{\textwidth}{ccccccc}
\midrule
\textbf{City} & \textbf{Radar} & \textbf{ADS-B} & \textbf{RF} & \textbf{Remote ID} & \textbf{Acoustic} & \textbf{Optical Camera} \\ [0.5ex]
\midrule
Columbus & 610 & 1 & 49 & 49 & 29500 & 55000 \\
Cleveland & 261 & 1 & 24 & 24 & 10684 & 21642 \\
Cincinnati & 228 & 1 & 20 & 20 & 20000 & 34335 \\
Akron & 153 & 1 & 14 & 14 & 13920 & 27228 \\
Toledo & 192 & 1 & 17 & 17 & 15305 & 33594 \\
Dayton & 110 & 1 & 11 & 11 & 10095 & 22890 \\
\midrule
\end{tabularx}
\end{table}
\vspace{-12pt}

% ADS-B, remote ID, and RF sensor types generate positive NPVs within the analysis period, as shown in Figure \ref{fig10}. ADS-B and Remote ID generate the largest NPV most quickly, while RF brings the third largest NPV. An increase in the range of a sensor will result in a larger coverage area. A higher generation of traffic surveillance data, requiring fewer sensors to cover a city. This, in turn, will raise the ratio of the coverage area to the cost of a sensor and ultimately lead to an increase in NPV. The projected NPV for the ADS-B sensor type steadily increases over time due to revenue generated from subscription fees and reach around \$18.50 million in the final year of the analysis horizon, as illustrated in Figure \ref{fig10a}. The projected NPV trend of the remote ID sensor type shows that it starts to generate positive NPV from 2024 and reaches approximately \$18.38 million in 2033. On the other hand, as shown in Figure \ref{fig10b}, the NPV for the RF sensor type is not positive until it reaches the BEP in 2029 and climbs to approximately \$13.91 million in 2033. The projected yearly NPVs for the radar, acoustic, and optical camera sensor types reveal negative NPVs due to the high initial capital cost for installing the sensors.

\begin{figure}[H]
{\captionsetup{position=bottom,justification=centering}
\begin{subfigure}{.48\textwidth}
  \centering
  % include first image
  \includegraphics[width=1\linewidth]{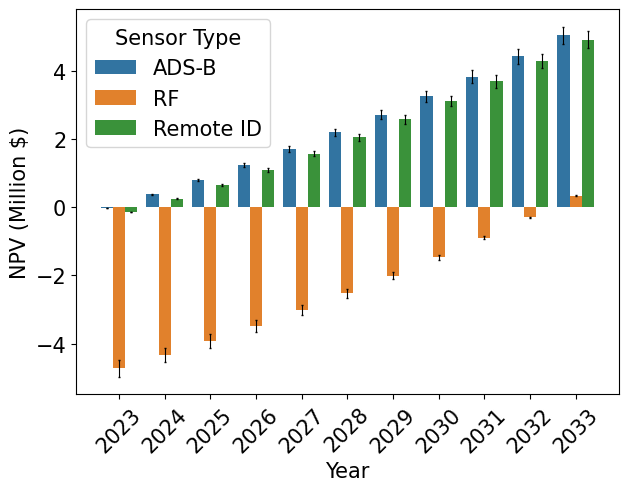}
  \centering
  \caption{Yearly NPV of ADS-B, RF, and remote ID.}
  \label{fig:NPVa}
\end{subfigure}
\begin{subfigure}{.48\textwidth}
  \centering
  % include second image
  \includegraphics[width=1\linewidth]{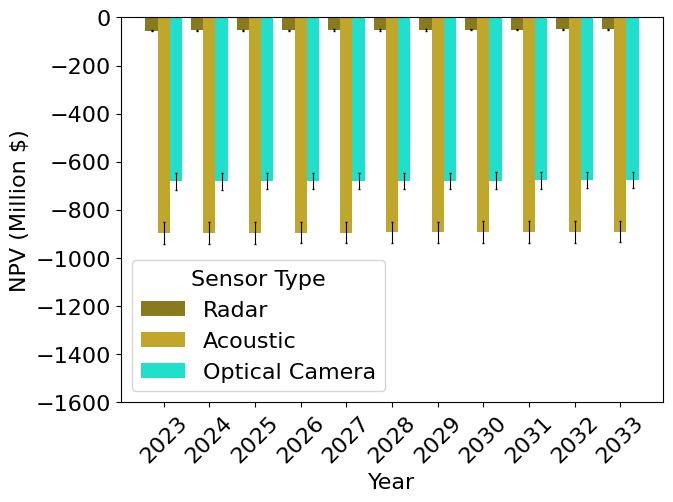}
  \centering
  \caption{Yearly NPV of radar, acoustic, and optical camera.}
  \label{fig:NPVb}
\end{subfigure}}
\caption{Yearly NPV of six sensor types.}
\label{fig:NPV}
\end{figure}

Utilizing the total sensor cost for SMCO presented in this section and the revenue and cloud computing cost for SMCO detailed in Section \ref{revenue_cloud}, the NPV is calculated. The yearly NPVs associated with all sensor types are presented in Figure \ref{fig:NPV}, where the length of an error bar for a sensor type in a given year, extending on either side of its central NPV value, indicates the range within which the NPV for that sensor type in that year is expected to fall. For all sensor types, the steady increase in NPV with time is fueled by the increase in yearly revenues generated from the subscription fees. This NPV growth is less noticeable for radar, optical camera and acoustic sensor types as they have high initial sensor costs. ADS-B, remote ID, and RF sensor types generate positive NPVs within the analysis period. ADS-B generates the largest NPV, followed closely by remote ID, while RF brings the third largest NPV. These sensor types lead the NPV race because they have lower unit prices and higher ranges, thus requiring fewer sensors to cover a city, and hence have lower sensor costs. Both the ADS-B and remote ID sensor types quickly reach BEP in 2024. Their projected NPVs reach around \$5.04 million and \$4.90 million in the final year of the analysis period, as illustrated in Figure \ref{fig:NPVa}. The RF sensor type takes longer to reach BEP, gaining positive NPV of approximately \$3.29 million in 2033. On the other hand, as shown in Figure \ref{fig:NPVb}, the projected yearly NPVs for the radar, acoustic, and optical camera sensor types feature negative NPVs over the 10-year analysis period due to their high initial sensor costs. 

\textbf{Policy recommendation}: As discussed previously in Section \ref{sensortype}, each sensor type can detect and track cooperative and/or non-cooperative aircraft flying in low-altitude airspace. Radar, RF, acoustic, and optical camera are capable of tracking both types of aircraft, while ADS-B and remote ID can only track cooperative aircraft. Based on the NPV, if tracking only cooperative aircraft is sufficient, then ADS-B and remote ID sensor types are recommended as they are the most financially viable sensor types for AAM surveillance network and LASIC. If tracking both cooperative and non-cooperative aircraft, especially those flying over penitentiary and other restricted areas, is a requirement, then RF is the most profitable sensor type.

\subsection{Heterogeneous Sensor Placement Analysis} \label{hetero}

% \begin{table}[H]
% \centering
% \caption{Comparison of total sensor cost between homogeneous and heterogeneous sensor networks.} \label{compare2}
% \label{compare2}
% \begin{tabular}{lccc}
% \midrule
% & \multicolumn{3}{c}{\textbf{Total Sensor Cost (million \$)}} \\\midrule
%  & \multicolumn{2}{c}{\textbf{Homogeneous Sensor Network}} & \textbf{Heterogeneous Sensor Network}\\
% \cmidrule{2-4}
% \textbf{City} & \textbf{Radar} & \textbf{Acoustic}  & \textbf{Radar and Acoustic} \\ \midrule
% Columbus & 21.35 & 265.50  & 20.75 \\
% Cleveland & 9.14 & 96.16  & 8.51 \\
% Cincinnati & 7.98 & 180.00  & 7.83 \\
% Akron & 5.36 & 125.28 & 5.06 \\
% Toledo & 6.72 & 137.75 & 6.72 \\
% Dayton & 3.85 & 90.86  & 3.52 \\ \midrule
% \end{tabular}
% \end{table}

%updated paragraph
The heterogeneous sensor placement analysis aims to investigate the network composition and costs associated with using a combination of sensors of different types rather than selecting sensors of just one type. The SAND model identifies the optimal sensor locations of the assorted sensor types to build the AAM surveillance sensor network across SMCO. To conduct an experiment on the heterogeneous sensor placement analysis, we consider providing coverage to sensitive locations within SMCO, such as penitentiaries, police stations, and airports, where detecting both cooperative and non-cooperative aircraft is equally important. Since radar, RF, acoustic, and optical camera sensor types can detect both types of aircraft, we initially consider the set of sensors as $S$ = \{Radar, RF, Acoustic, Optical Camera\}. However, we find that only the RF sensor type is selected from this set, as it dominates other sensor types and generates the same results as the RF sensor type in the homogeneous case (refer to Section \ref{homo}). This is due to RF sensor type having a larger range, higher detection probability on all terrain types, higher field of view (lower $\delta$), and a lower unit price compared to other sensor types in the set. We then consider the set of sensors as $S$ = \{Radar, Acoustic, Optical Camera\} to conduct the experiment again and generate further insights on the heterogeneous sensor placement analysis, which are given in this section.

%old paragraph
% The heterogeneous sensor placement analysis aims to investigate the network composition and costs associated with using a combination of sensors of different types rather than selecting sensors of just one type. The SAND model identifies the optimal sensor locations of the assorted sensor types to build the AAM surveillance sensor network across SMCO. This analysis is particularly useful when it comes to providing coverage to the sensitive locations within SMCO, such as penitentiaries, police stations, and airports, where detecting both cooperative and non-cooperative aircraft are equally important. For this analysis, radar, acoustic, and optical camera sensor types were considered as they can detect both types of aircraft. Also, the analysis focused on only the city of Akron. 

The total number of sensors and sensor cost needed to place the sensors, as presented in Table \ref{compare1} and Table \ref{compare2}, respectively, are compared between the two types of sensor placement: homogeneous sensor network and heterogeneous sensor network. For each city, the homogeneous cases show the values of total number of sensors and sensor cost of placing three individual sensor types (radar, acoustic, and optical camera), whereas the heterogeneous case shows the values for a mix of these sensor types. Comparing the values shows that heterogeneous sensor placement requires fewer sensors than acoustic sensors and optical cameras in a homogeneous sensor network but more sensors than radars in a homogeneous sensor network. For example, in Akron, the heterogeneous sensor network requires 155 sensors as shown in Table \ref{compare1}, which is lower than the numbers required for acoustic and optical cameras, at 13920 and 27228 respectively, but higher than the number of sensors required for radars in the homogeneous sensor network, which is 153. However, the total cost of sensors for setting up a heterogeneous sensor network is much lower than the separate costs for radar, optical camera, and acoustic sensor types in the homogeneous case, which are \$1.89 million, \$13.53 million, and \$31.56 million, respectively. The SAND model selects the optimal number and location of sensors to minimize the total cost of sensors, even if it means using more sensors by replacing some radars with acoustic sensors. Although this results in an increased number of sensors, the use of lower-priced sensors reduces the total sensor cost, which is the objective of the SAND model. 

\begin{table}[H]
\caption{Comparison of total number of sensors needed between homogeneous and heterogeneous sensor networks.} \label{compare1}
\newcolumntype{C}{>{\centering\arraybackslash}X}
\setlength{\tabcolsep}{2.75mm}
\begin{tabularx}{\textwidth}{lcccc}
\midrule
& \multicolumn{4}{c}{\textbf{Total Number of Sensors Needed}} \\\midrule
 & \multicolumn{3}{c}{\textbf{Homogeneous Sensor Network}} & \textbf{Heterogeneous Sensor Network}\\
\cmidrule{2-5}
\textbf{City} & \textbf{Radar} & \textbf{Acoustic} & \textbf{Optical Camera} & \textbf{Radar, Acoustic, and Optical Camera} \\ \midrule
Columbus & 610 & 29500 & 55000 & 648 \\
Cleveland & 261 & 10684 & 21642 & 333 \\
Cincinnati & 228  & 20000 & 34335 & 240 \\
Akron & 153 & 13920 & 27228 & 155 \\
Toledo & 192 & 15305 & 33594 & 192 \\
Dayton & 110 & 10095 & 22890 & 112 \\ \midrule
\end{tabularx}
\end{table}
\vspace{-12pt}

\begin{table}[H]
\caption{Comparison of total sensor cost between homogeneous and heterogeneous sensor networks.} \label{compare2}
\newcolumntype{C}{>{\centering\arraybackslash}X}
\setlength{\tabcolsep}{2.75mm}
\begin{tabularx}{\textwidth}{lcccc}
\midrule
& \multicolumn{4}{c}{\textbf{Total Sensor Cost (million \$)}} \\\midrule
 & \multicolumn{3}{c}{\textbf{Homogeneous Sensor Network}} & \textbf{Heterogeneous Sensor Network}\\
\cmidrule{2-5}
\textbf{City} & \textbf{Radar} & \textbf{Acoustic} & \textbf{Optical Camera} & \textbf{Radar, Acoustic, and Optical Camera} \\ \midrule
Columbus & 21.35 & 265.50 & 192.50 & 20.75 \\
Cleveland & 9.14 & 96.16 & 75.75 & 8.51 \\
Cincinnati & 7.98 & 180.00 & 120.17 & 7.83 \\
Akron & 5.36 & 125.28 & 95.30 & 5.06 \\
Toledo & 6.72 & 137.75 & 117.58 & 6.72 \\
Dayton & 3.85 & 90.86 & 80.12 & 3.52 \\ \midrule
\end{tabularx}
\end{table}

% Reductions in the number of sensors and cost vary across cities, with Akron showing the highest reduction and Dayton the lowest. 

% Therefore, deploying a mix of sensors (Heterogeneous) can reduce the cost of sensor deployment compared to deploying individual sensor types (Homogeneous).

% \begin{figure}[H]
%     \centering
%     \includegraphics[width=10cm,height=10cm]{akron_mix2.png}
%     \caption{Optimal locations of mixed sensors in Akron (red markers represent the location of radars, the green markers the location of acoustic sensors, and the texts below the markers show the number of sensors needed at the respective locations).}
%     \label{fig:mix}
% \end{figure}

\begin{figure}[H]
{\captionsetup{position=bottom,justification=centering}
\begin{subfigure}{.49\textwidth}
  \centering
  % include first image
  \includegraphics[width=1\linewidth]{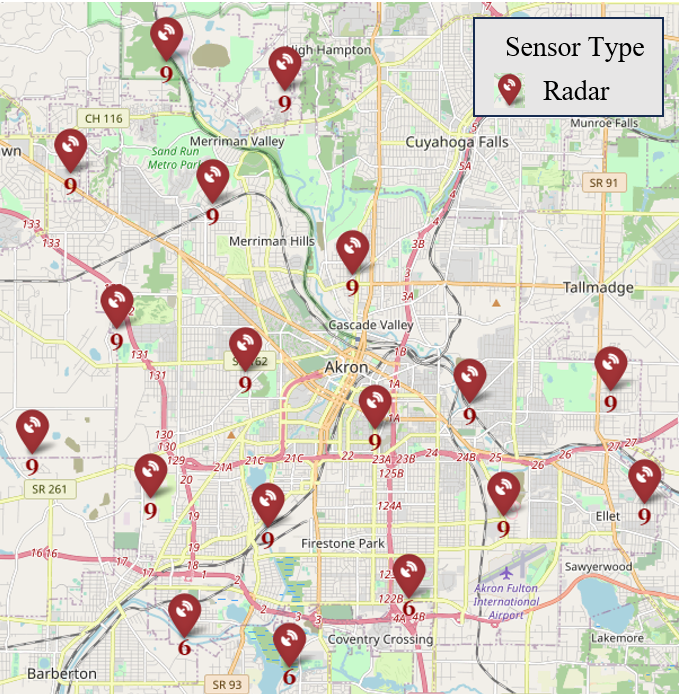}
  \centering
  \caption{Optimal locations of sensors in homogeneous network.}
  \label{fig:mixa}
\end{subfigure} \hfill
\begin{subfigure}{.49\textwidth}
  \centering
  % include second image
  \includegraphics[width=1\linewidth]{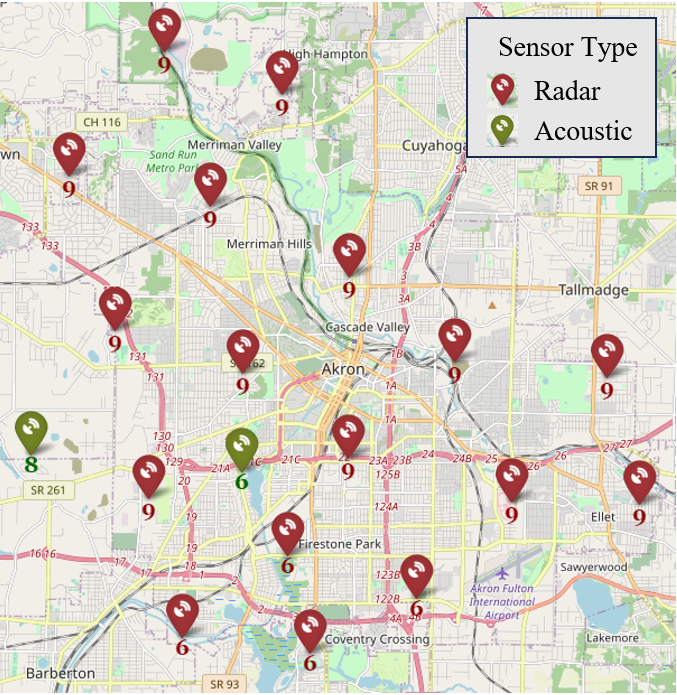}
  \centering
  \caption{Optimal locations of mixed sensors in heterogeneous network.}
  \label{fig:mixb}
\end{subfigure}}
\caption{Optimal locations of sensors in Akron in homogeneous and heterogeneous sensor networks (texts below the markers show the number of sensors needed at the respective locations).}
\label{fig:mix}
\end{figure}

The SAND model reconfigures the number and location of sensors in the heterogeneous sensor network, as demonstrated in Figure \ref{fig:mix}, to optimize the total sensor cost more effectively. The figure shows the optimal sensor locations in both the homogeneous and heterogeneous sensor networks for the city of Akron, where red markers represent radar locations, green markers represent acoustic sensor locations, and the numbers below the markers indicate the number of sensors needed at each location. In the heterogeneous sensor network, an optimal combination of sensors with appropriate ranges is chosen to cover the given area, resulting in lower costs compared to the homogeneous sensor network. In other words, the ranges of the different sensor types are utilized effectively in the heterogeneous sensor network to reduce the sensor cost. For example, near the outer edges of the city and in small pockets within the city, sensors with a smaller range and lower unit price are placed, such as acoustic sensors, instead of sensors with a higher range and higher unit price, like radar, to minimize the sensor cost.

\begin{figure}[H]
    \includegraphics[width=9cm,height=7cm]{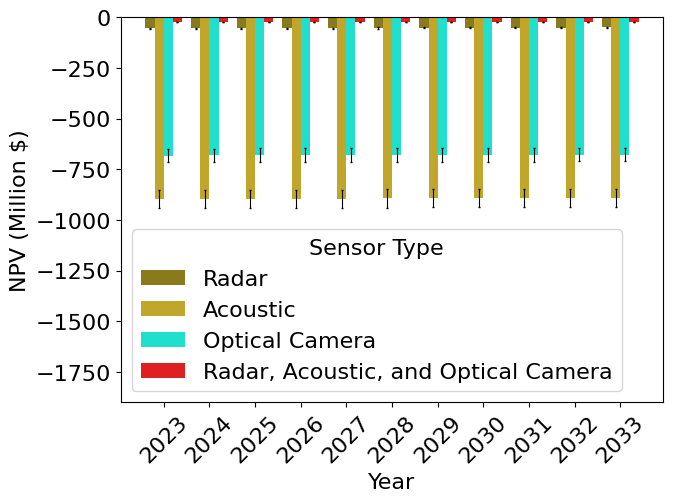}
    \caption{Comparison of yearly NPV between the homogeneous and heterogeneous sensor networks.}
    \label{fig:hetero}
\end{figure}

We then calculate NPV of the heterogeneous sensor network in SMCO by taking into account the revenue, cloud computing cost, and total sensor cost. As shown in Figure \ref{fig:hetero}, the NPV generated by the heterogeneous sensor network is higher than the NPVs generated by the respective homogeneous sensor networks for individual radar, acoustic sensor, and optical camera. This is due to the lower sensor cost of the heterogeneous sensor network, which is discussed earlier in this section. The first three bars in the figure represent the homogeneous cases, while the last one represents the heterogeneous case with a mix of radar, acoustic sensor, and optical camera.Although the NPVs for the radar, acoustic, and optical camera mentioned in Section \ref{homo} are negative, there is a noticeable trend of increasing NPV values as we move from the homogeneous cases to the heterogeneous case in all time periods. Therefore, we conclude that a heterogeneous sensor placement analysis can help generate more NPV than a homogeneous sensor placement analysis.

\textbf{Policy recommendation}: The heterogeneous sensor network offers a lower sensor cost compared to a homogeneous sensor network, allowing for the identification of the optimal mix of sensors from a given set of sensor types. Therefore, among the two types of sensor placement, a heterogeneous sensor network is recommended when the right set of sensor types is selected, considering which types of aircraft need to be detected and tracked. By using heterogeneous sensor placement, it is possible to design a surveillance network with minimum cost which installs sensor types that can track non-cooperative aircraft (e.g., RF) in security sensitive areas (e.g., penitentiaries, law enforcement facilities, and correctional facilities) and sensor types that can track either cooperative or non-cooperative aircraft (e.g., ADS-B and remote ID) in non-security sensitive or general public areas. Other sensor placement constraints can also be enforced while designing heterogeneous sensor network based on the requirements, preferences and regulations of the government, AAM operators, and LASIC subscribers. 

% In the heterogeneous scenario, a lower number of sensors was needed as compared to the homogeneous sensor analysis. The results show that heterogeneous helps to minimize the number of sensors needed and, hence, the total sensor cost, by choosing an optimal combination of sensors within an area. Since radar provides a larger coverage area compared to the other two sensor types, fewer radars are needed to cover a city. However, to minimize the cost and cover the blocks near the outer edges of the city, the model tries to place sensors with a smaller range, such as optical cameras, instead of sensors with a higher range, like radar. The cost of an individual acoustic sensor is inexpensive, but it requires a large quantity to cover a specific area as it has a smaller range compared to an optical camera. Hence, optical camera is able to minimize the the sensor cost more effectively, resulting in their selection in the optimal mix of heterogeneous sensors.

\subsection{Sensitivity Analysis}

% A sensitivity analysis was conducted to evaluate the impact of changes in subscription fees and the initial number of potential subscribers on NPV over the time horizon. Based on survey responses, three values for the monthly subscription fee per subscriber ($S$), \$100, \$250, and \$400, and three values for the number of potential subscribers in 2024 ($N$), 50, 75, and 100 were considered.

% The trend of NPV, as illustrated in Figure \ref{fig11a}, shows that an increase in the subscription fee leads to a corresponding increase in the projected yearly NPV. The NPV in 2033 ranges from \$4 million to \$17 million when the subscription fee per subscriber ($S$) varies from \$100 to \$400. This change also affects the BEP, as demonstrated by the example of RF sensors. The BEP for RF sensors changes from 2029 to 2031 when the subscription fee changes from \$400 to \$250.

% Similarly, as depicted in Figure \ref{fig11b}, the projected yearly NPV increases with the increase in the number of subscribers. The NPV in 2033 ranges from \$8 million to \$17 million when the number of subscribers ($N$) varies from 50 to 100, and the BEP for RF sensors changes from 2033 to 2031 accordingly.

Despite the growing interest in AAM, several unresolved obstacles and concerns exist, including the need for new widespread infrastructure to support AAM operations, such as vertiports, takeoff and landing sites, charging stations, air traffic control systems, airspace routes, and surveillance networks. Additionally, the current regulatory framework for air transportation is not designed for AAM, requiring new regulations and standards. Factors such as changes in consumer preferences, regulatory requirements, and technological advancements could all affect the adoption and growth of AAM services. Hence, in this section, a sensitivity analysis is conducted to examine the impact of changes in parameters directly affecting NPV generated in our analysis, such as LASIC subscription fee and the number of initial subscribers. We also conduct sensitivity analysis for other key parameters that we are uncertain about due to the lack of verified data available, such as terrain-based sensor detection probabilities and minimum required detection probability, to observe their impacts on respective outputs.

\begin{figure}[H]
\begin{subfigure}{.48\textwidth}
  \centering
  % include first image
  \includegraphics[width=1\linewidth]{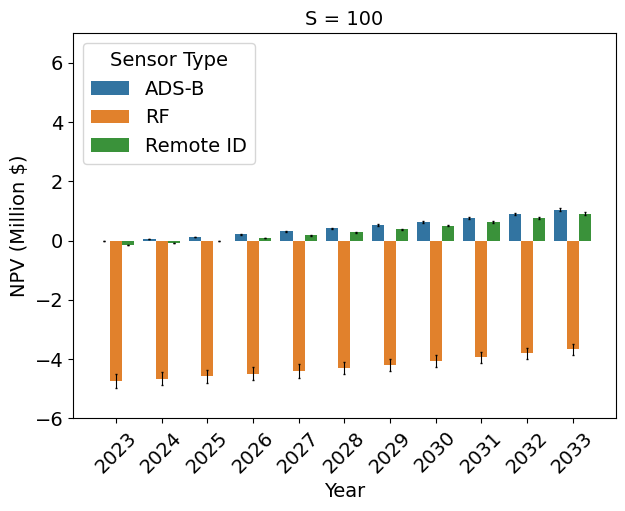}
\end{subfigure}
\begin{subfigure}{.48\textwidth}
  \centering
  % include second image
  \includegraphics[width=1\linewidth]{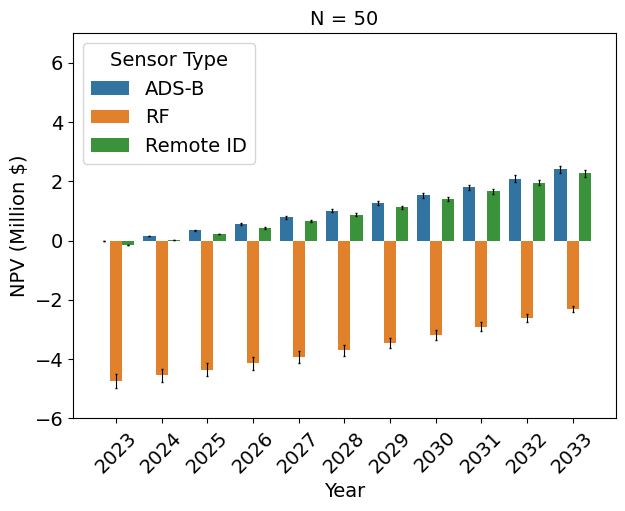}  
\end{subfigure}
\newline
\begin{subfigure}{.48\textwidth}
  \centering
  % include third image
  \includegraphics[width=1\linewidth]{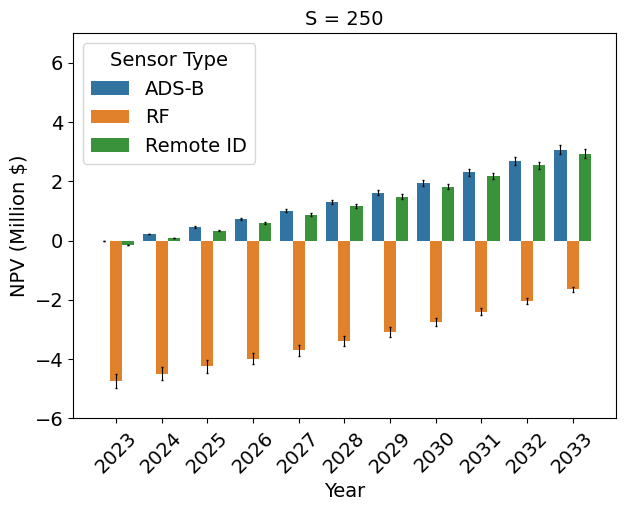}  
\end{subfigure}
\begin{subfigure}{.48\textwidth}
  \centering
  % include fourth image
  \includegraphics[width=1\linewidth]{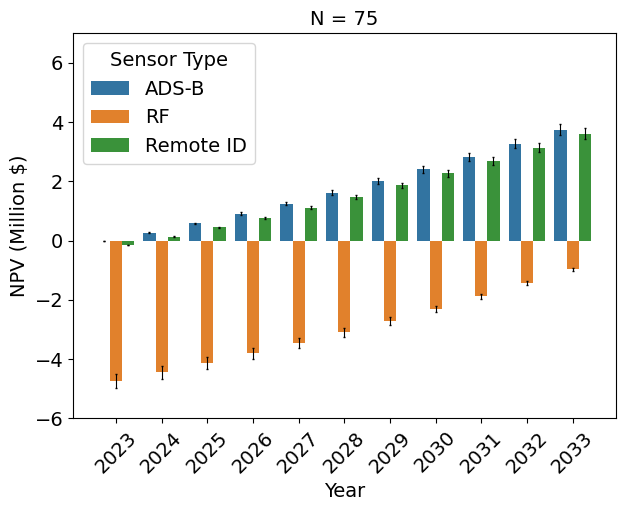}  
\end{subfigure}
\newline
\begin{subfigure}{.48\textwidth}
  \centering
  \includegraphics[width=1\linewidth]{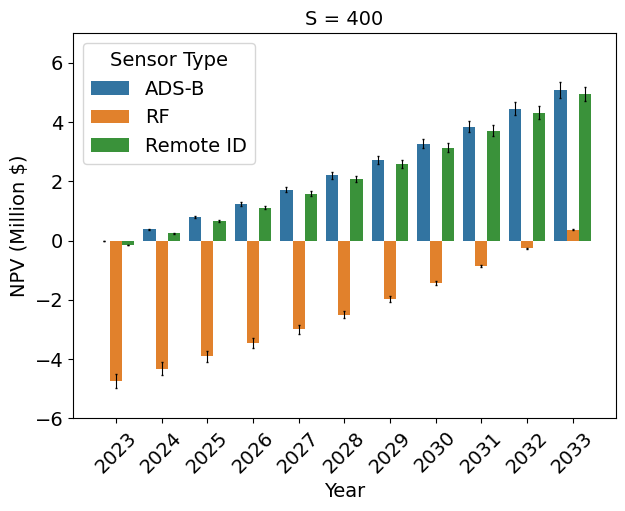}  
  \caption{Yearly NPV for six sensor types varying $S$.}
  \label{fig11a}
\end{subfigure}
\begin{subfigure}{.48\textwidth}
  \centering
  \includegraphics[width=1\linewidth]{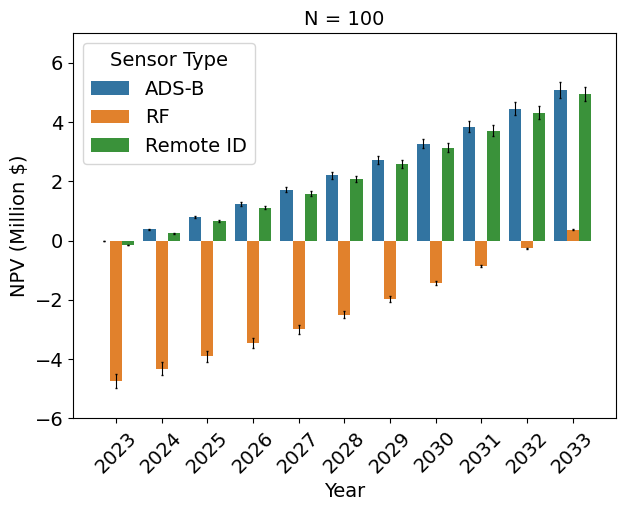}  
  \caption{Yearly NPV for six sensor types varying $N$.}
  \label{fig11b}
\end{subfigure}
\caption{Yearly NPV for six sensor types varying $S$ and $N$.}
\label{fig:sensitivity}
\end{figure}

To address the uncertainty associated with the demand for AAM, it is crucial to consider a range of subscription fees and yearly numbers of subscribers, instead of fixed values, as these variables directly affect the revenue and NPV generated by LASIC. Therefore, a sensitivity analysis is conducted to examine the impact of changes in these variables on the NPV over the analysis period. For this analysis, the sensor types with the three highest producing NPVs in homogeneous sensor network --- ADS-B, RF and remote ID --- are considered. To vary the yearly number of subscribers for LASIC, the number of subscribers in the initial year (2024) is varied, which affects the number of subscribers in the subsequent years. Based on survey responses, three values for the monthly subscription fee per subscriber ($S$) --- \$100, \$250, and \$400 --- and three values for the number of potential subscribers in 2024 ($N$) --- 50, 75, and 100 --- are considered. The trends observed in Figure \ref{fig:sensitivity} show that higher values of $S$ and $N$ lead to increases in the NPV and causes the BEP to occur earlier. These effects are attributed to the increase in revenue generation prompted by increases in $S$ and $N$. As demonstrated by the example of RF sensor type, when $S$ is \$100, the NPV in 2033 is \$-3.67 million and cannot reach BEP within the analysis period. On the other hand, when $S$ increases to \$400, for the same sensor type, NPV reaches BEP of \$0.37 million by 2033. Similarly, if $N$ is 50, the NPV shows a net loss of \$2.31 million in 2033, and it is not possible to reach BEP during the analysis period. However, for the same sensor type, when $N$ increases to 100, the NPV is expected to reach BEP of \$3.74 million by 2033. The analysis highlights that the state government can achieve a net profit for several sensor types within the analysis period as long as $S$ and $N$ are not too low. The values of $S$ and $N$ at which positive NPV is ensured are contingent on the chosen sensor type. If the objective is to detect and track solely cooperative aircraft, then ADS-B and remote ID are profitable options as the sensor type for the AAM surveillance network and LASIC, as discussed in Section \ref{homo}. For this case, even if $S$ and $N$ assume values of \$100 and 50, respectively, the state government can still attain a net profit within the analysis period. On the other hand, if both cooperative and non-cooperative aircraft are required to be tracked, then RF represents the most profitable sensor type for the AAM surveillance network and LASIC, as discussed in Section \ref{homo}. In this scenario, the state government should ensure $S$ to be no less than \$400 when $N$ is 100 to achieve a net profit within the analysis period.

To understand the impact of changes in terrain-based sensor detection probabilities on the number of required sensors and the total sensor cost, it is important to consider that sensor detection probabilities can vary depending on the type and manufacturer of the sensor. For instance, if we focus on radar sensor type and consider the city of Akron, we can analyze the effect of varying the sensor detection probabilities for different terrain types by conducting a sensitivity analysis. The sensor detection probabilities for different terrain types in Akron and the corresponding number of sensors and sensor cost are presented in Table \ref{radar:sen}, along with the effects of a 5\% increase and 5\% decrease in sensor detection probabilities. The preset value of sensor detection probabilities ranges from 0.75 to 0.95 for different terrain types (refer to Table \ref{tab2}). A 5\% increase in sensor detection probabilities results in a decrease in both the number of sensors required and total sensor cost. This is because higher sensor detection probability allows the radar to cover the same area with fewer sensors while increasing the probability of tracking AAM aircraft, thereby reducing the sensor cost for a given sensor type. Conversely, a 5\% decrease in sensor detection probabilities results in an increase in both the number of sensors required and total sensor cost, as a lower sensor detection probability requires more sensors to cover the same area.

The sensitivity analysis of the minimum required detection probability is important because it helps to understand how the number of sensors and their cost vary with changes in $r$. A higher $r$ value means that the sensor needs to be more sensitive to detect AAM aircraft with a higher level of confidence, which can increase the number of sensors required and the cost of the sensors. For example, in Table \ref{rt} provided for Akron City, we can see how the number of sensors and their cost vary for different sensor types as $r$ increases. The number of sensors and the cost of a sensor type are more affected when the sensor type has a lower range, lower detection probability on all terrain types, lower field of view (higher $\delta$), and a higher unit price. For example, the acoustic sensor type is the most affected by changes in $r$. As the value of $r$ increases from 0.96 to 0.99, the cost of the acoustic sensor type increases significantly from \$108.612 million to \$154.350 million and the number of sensors required also increases from 12,068 to 17,150. On the other hand, the ADS-B sensor type is the least affected by changes in $r$. We can see that for ADS-B sensors, increasing the required detection probability from 0.96 to 0.99 leads to an increase in the cost of the sensors from \$0.002 million to \$0.005 million, and the number of sensors required from 1 to 2. 

\vspace{-6pt}

\begin{table}[H]
\caption{\label{radar:sen} Detection probability of a radar for different terrain types in Akron and the corresponding number of sensors and sensor cost, along with the effects of a 5\% increase and 5\% decrease in detection probability.}
\begin{adjustwidth}{-\extralength}{0cm}
		\newcolumntype{C}{>{\centering\arraybackslash}X}
		\setlength{\tabcolsep}{3.5mm}
		\begin{tabularx}{\fulllength}{lccccccr}
\midrule
&&&Terrain Type ($T$)& \multicolumn{2}{c}{} & \\ \cmidrule{2-6}
& Open & Water & Neighborhood & Hill & Commercial & Number of & Total Sensor Cost \\
Case & & & & & Area &  Sensors & (\$ million) \\\midrule
Preset Value & 0.9500 & 0.9000 & 0.8500 & 0.7500 & 0.7500 & 153 & 5.355\\
5\% Increase & 0.9975 & 0.945 & 0.8925 & 0.7875 & 0.7875 & 114 & 3.990 \\
5\% Decrease & 0.9025 & 0.855 & 0.8075 & 0.7125 & 0.7125 & 165 & 5.775\\ \midrule
\end{tabularx}
\end{adjustwidth}
\end{table}

\begin{table}[H]
\caption{Number of sensors and cost for different sensor types at different $r$ values.}
\label{rt}
\newcolumntype{C}{>{\centering\arraybackslash}X}
\setlength{\tabcolsep}{3.5mm}
\begin{tabularx}{\textwidth}{ccccc}
\midrule
\textbf{Sensor Type} & $\boldsymbol{r}$ & \textbf{Number of Sensors} & \textbf{Total Sensor Cost (million \$)} \\ [0.5ex] 
\midrule
ADS-B & 0.96 & 1 & 0.002 \\ 
& 0.97 & 1 & 0.002 \\
& 0.98 & 1 & 0.002 \\
& 0.99 & 2 & 0.005 \\ 
\midrule
Remote ID & 0.96 & 10 & 0.011 \\ 
& 0.97 & 11 & 0.012 \\
& 0.98 & 14 & 0.015 \\
& 0.99 & 15 & 0.017 \\ 
\midrule
Radar & 0.96 & 117 & 4.095 \\ 
& 0.97 & 120 & 4.200 \\
& 0.98 & 153 & 5.355 \\
& 0.99 & 165 & 5.775 \\ 
\midrule
RF & 0.96 & 10 & 0.350 \\ 
& 0.97 & 11 & 0.385 \\
& 0.98 & 14 & 0.490 \\
& 0.99 & 15 & 0.525 \\ 
\midrule
Acoustic & 0.96 & 12068 & 108.612 \\ 
& 0.97 & 12421 & 111.789 \\
& 0.98 & 13920 & 125.280 \\
& 0.99 & 17150 & 154.350 \\ 
\midrule
Optical Camera & 0.96 & 19986 & 69.951 \\ 
& 0.97 & 26742 & 93.597 \\
& 0.98 & 27228 & 95.298 \\
& 0.99 & 28680 & 100.380 \\ 
\midrule
\end{tabularx}
\end{table}

\section{Conclusion and Future Work}\label{sec5}

To enable real-time detection and tracking of AAM aircraft flying at lower altitudes, an effective AAM surveillance network is required to ensure adequate coverage and monitoring. This study introduces the SAND model, which aims to design an AAM surveillance network that provides full coverage in a specified operational area while minimizing the total sensor cost. The model considers various factors such as sensor types, terrain types, terrain-based sensor detection probabilities, and minimum detection probability requirement. We consider two types of surveillance sensor networks: homogeneous sensor network consisting of only one type of sensor, and heterogeneous sensor network consisting of different types of sensors. Additionally, we present LASIC as a centralized cloud database that needs to be connected to the AAM surveillance network to efficiently store and process the large amounts of data generated by the network. The required features and functionalities of LASIC are determined based on AAM market data and survey inputs. To justify the investment in AAM surveillance network and LASIC, a rigorous data-driven cost-benefit analysis is conducted by identifying, quantifying, and evaluating the costs and benefits associated with the infrastructure. We conducted the analysis for the State of Ohio over a 10-year period by estimating NPV for different sensor types. 

% To address the challenge of managing and utilizing the considerable number of data generated by AAM traffic, Then LASIC is presented %It aims  to achieve improved coordination among users of low altitude airspace.

% The analysis identifies two major cost factors: surveillance sensor cost and cloud computing cost, while the benefit factor is the revenue generated, which depends on the subscription fee and number of potential subscribers of LASIC.

The cost-benefit analysis identifies two significant cost factors --- surveillance sensor cost and cloud computing cost --- along with a benefit factor, which is the revenue generated by LASIC. This revenue is influenced by the subscription fee and the number of potential subscribers. The cloud computing cost for LASIC depend on the cloud server pricing policies and the required components for real-time and offline features. These costs are influenced by the projected AAM traffic and resulting surveillance data generated in each city. Due to larger areas and higher air traffic demand forecasts in Cleveland, Columbus, and Cincinnati compared to Akron, Toledo, and Dayton, the cloud computing costs are higher in the former set of cities. According to the homogeneous sensor placement analysis, the most profitable sensor types for detecting cooperative aircraft are ADS-B and remote ID sensor types, whereas for tracking both cooperative and non-cooperative aircraft, the most profitable option is the RF sensor type. This is because these sensors have a larger range, higher field of view, higher detection probability, lower unit price, and lower cost, leading to higher NPVs compared to other sensor types. The findings also indicate that by selecting an optimal combination of sensors of different sensor types to effectively cover a given area, a lower cost and higher NPV can be achieved through heterogeneous sensor placement. Furthermore, the results show that the total sensor cost in each city of SMCO varies based on factors such as city area, shape, terrain type, and terrain-based sensor detection probabilities. Given these factors, Columbus was found to require the largest sensor cost, while Dayton has the lowest sensor cost.

% In this study,  

% The SAND model determined which of the available locations within each city of SMCO should be selected to place sensors in different terrains of SMCO and track the aircraft with a minimum required detection probability. 

Because of the uncertainties in AAM demand and the significant influence of certain parameters on the results, such as the LASIC subscription fee, number of initial subscribers, terrain-based sensor detection probabilities, and minimum required detection probability, we perform a sensitivity analysis. This analysis aims to observe how changes in these parameters impact the results. The analysis indicates that an increase in terrain-based sensor detection probabilities leads to a decrease in the required number of sensors and total sensor cost, while a decrease in detection probability has the opposite effect. The analysis also reveals that an increase in minimum required detection probability leads to an increase in both the number of sensors required and total sensor cost, with the acoustic sensor type being most affected and the ADS-B sensor type being least affected. Furthermore, the analysis show that higher values of subscription fees and numbers of subscribers lead to increases in the NPV generated by LASIC and cause the BEP with respect to NPV to occur earlier, as they increase revenue generation. 

% , where the focus can be to cover maximum area based on the terrain features of each corridor
% Lastly, our another potential future work can be focused on developing and validating a risk assessment model to investigate the safety benefits of LASIC in AAM infrastructure. To achieve this, a simulation of various AAM operation scenarios, with and without LASIC, can be performed, and the resulting collision risk estimates can be compared. Such an analysis would provide valuable insight into the impact of implementing LASIC on AAM safety. \hl{In future work, the assumptions made during this study can be relaxed. For instance, the consideration of sensor obstruction by tall buildings could be reconsidered.

This study has produced several insights related to AAM surveillance network and LASIC and opportunities for future research that we plan to explore further. In future work, it is crucial to consider relaxing the assumptions made during this study, as doing so would render the solutions generated by the SAND model more practical. For instance, reevaluating the consideration of potential sensor obstructions will help ensure that detection requirements are met, even in the presence of natural or human-made structures in the AAM surveillance area. Furthermore, it is essential to account for trajectory planning, as the trajectories of AAM flights determines the density of flights in the airspace. Sensors have a certain capacity of aircraft that they can detect at any given time within their range. The sensor placement solution should factor this in to avoid exceeding their capacities. Within this trajectory planning, the SAND model can strategically position sensors based on the expected distances between the trajectories and sensor locations, taking into account that sensor detection performance depends on the distance between the sensor and the aircraft. Furthermore, other factors that affect sensor detection performance, such as weather conditions and sensor failures, can be considered in future work. Additional potential future work includes situations where there is a need to ensure surveillance with high probability of detection over a large area while working with limited resources, a different objective function can be considered to set up a surveillance network that maximizes the detection probability under a fixed budget for the sensors. This approach would be particularly important for border security, disaster response, or wildlife monitoring, where maintaining a high level of vigilance is crucial, but the cost of setting up and maintaining the surveillance network must be balanced with limited budget. Moreover, this model can be applied for designing surveillance networks for conventional air traffic as well. Next, the location of the sensors along and across the AAM corridors connecting the major cities of a state can be identified using the model. Furthermore, the optimal sensor placement model and LASIC framework presented in this paper can be applied to other potential states of the USA where AAM is emerging and has the potential for AAM demand. This will facilitate the development of an optimal surveillance network and surveillance information clearinghouse for AAM in these states. However, it is essential to evaluate the impacts of such investments by conducting a cost-benefit analysis similar to the one performed in this study before making any significant investments in AAM surveillance infrastructure. Another potential future work can be focused on developing and validating a risk assessment model to investigate the safety benefits of LASIC in AAM infrastructure. To achieve this, a simulation of various AAM operation scenarios, with and without LASIC, can be performed, and the resulting collision risk estimates can be compared. Such an analysis would provide valuable insight into the impact of implementing LASIC on AAM safety.

\vspace{6pt} 

\acknowledgments{This research was supported by the Ohio Department of Transportation (Agreement \mbox{No.: 36496}, PID: 114242, SJN: 136337), and the system requirements study was led by CAL Analytics under the same project. The authors of the article acknowledge and express gratitude for the input of the many people directly involved and consulted during the performance of this effort. Special gratitude for guidance on this project goes to the Technical Advisory Committee members, Fred Judson and Richard Fox, and Sean Calhoun from CAL Analytics and Rubén Del Rosario from Crown Consulting.}

\conflictsofinterest{The authors declare no conflict of interest.} 

%%%%%%%%%%%%%%%%%%%%%%%%%%%%%%%%%%%%%%%%%%
\begin{adjustwidth}{-\extralength}{0cm}
%\printendnotes[custom] % Un-comment to print a list of endnotes

\reftitle{References}

\end{adjustwidth}

\end{document}